\documentclass[a4paper,10pt]{article}

\usepackage{amsmath,amsfonts,amscd,amssymb,amsthm}
\usepackage[all]{xy}
\usepackage[dvips]{graphics}

\newtheorem{thm}{Theorem}[section]
\newtheorem{lemme}[thm]{Lemma}
\newtheorem*{Thm}{Theorem}
\newtheorem{souslemme}[thm]{Sub-lemma}
\newtheorem{definition}{Definition}[section]
\newtheorem{coro}{Corollary}[thm]

\newtheorem{prop}{Proposition}[thm]
\newtheorem{propr}{Property}[thm]
\newcommand{\dd}{\text{d}}
\newcommand{\OO}{\mathcal{O}}
\newcommand{\DD}{\mathcal{D}}
\newcommand{\XX}{\mathcal{X}}
\newcommand{\FF}{\mathcal{F}}

\newcommand{\ZZ}{\mathcal{Z}}
\newcommand{\MM}{\mathcal{M}}
\newcommand{\UU}{\mathcal{U}}

\newcommand{\GG}{\mathcal{G}}
\newcommand{\NN}{\mathcal{N}}
\newcommand{\PP}{\mathcal{P}}
\newcommand{\fami}[1]{\left\{ #1 \right\}}
\newcommand{\coup}[1]{\left( #1 \right)}

\newcommand{\carre}{
\begin{flushright}
$\square$
\end{flushright}
\medskip}
\newcommand{\dbar}{\overline{\mathbb{D}}}
\newcommand{\fraction}[2]{\displaystyle\frac{#1}{#2}}
\newenvironment{demo}{\noindent {\it Proof :}}{\carre}
\newcommand{\fleche}[1]{\stackrel{\scriptscriptstyle #1}{ \xymatrix{\ar@{-->}[r]&  \\ }}}

\begin{document} 

\title{Construction of foliations with prescribed separatrix.}

\author{Y. Genzmer}
\maketitle
%

\begin{abstract} We build a germ of singular foliation in $\mathbb{C}^2$ with analytical class of separatrix and holonomy representations prescribed. Thanks to this construction, we study the link between moduli of a foliation and moduli of its separatrix.
\end{abstract}

\section*{Introduction}Considering the problem of moduli for a germ of singular holomorphic foliation $\FF$ in $\mathbb{C}^2$ leads naturally to point out many kinds of topological and analytical invariants. The invariants of first kind come from the reduction of singularities $E:(\MM,\DD)\rightarrow (\mathbb{C}^2,0)$ of the foliation:
\begin{enumerate}
\item the topological class of the manifold $\MM$; this invariant is a combinatorial one;
\item the analytical class of the pointed divisor $\DD$; 
\item the analytical class of the manifold $\MM$;
\end{enumerate}
The invariants of second kind are more related to the foliation itself: the collection of projective holonomy representations defined over each component of the divisor $\DD$, the so called holonomy pseudo-group.
 A natural problem is to know if \emph{coherent} data of above invariants correspond to a concrete foliation. The first step toward this question is a theorem of A. Lins-Neto \cite{alcides}, which establishes the possibility of constructing a foliation with invariant $(1)$ and projective holonomy prescribed. In his thesis \cite{Seguy}, M. Seguy shows that it is possible to fix the invariant $(2)$. The aim of this article is to prove that one can even prescribe the invariant $(3)$ and the holonomy invariants in the construction of a foliation. 

\noindent The three first sections of this paper are devoted to prove the above result. The main tool of our construction is the equisingular unfolding of foliation. Basically, it's an iso-holonomic topogically trivial deformation. It was introduced in \cite{univ} to study the local moduli of a foliation. It is not easy to find in a constructive way such a deformation not analytically trivial: except in the very special case where the foliation admits a first integral, none exemple is known. This leads one to study the cohomological interpretation in order to build non-trivial equisingular unfolding. In the first section, we fix a manifold $\MM$, which is built over the reduction tree of a foliation and we define a very special class of manifold denoted by $\text{Glu}_0(\MM,\UU,Z)$ related to $\MM$. In the second section, $\MM$ is supposed to be foliated by $\FF$. A property of cobordism type is pointed out and allows us to detect the existence of a foliation on any element of $\text{Glu}_0(\MM,\UU,Z)$. This foliation will automatically be linked to $\FF$ by an equisingular unfolding. The theorem (\ref{thmA}) states that, under the generical hypothesis of being of {\it second kind} (\ref{deuxspe}),  the cobordism property holds for any element of $\text{Glu}_0(\MM,\UU,Z)$. In the third section, we deduce the following result from the cobordism property
\begin{Thm}[Existence theorem]
\noindent Let $\FF$ be supposed of second kind. If $\MM$ and $\MM'$ are topologically equivalent then there exists an holomorphic foliation on $\MM'$ linked to $\FF$ by an equisingular unfolding.
\end{Thm}
\noindent A trivial but maybe more explicit corollary of the above theorem is the following:
\begin{Thm}
Let $\omega_0$ be a germ of singular holomorphic $1$-form  of second kind at $0\in\mathbb{C}^2$ and $E:\MM\rightarrow \mathbb{C}^2$ its singularities reduction. For any blowing-up process $E'$ with same dual tree as $E$, there exists a $1$-form $\omega'$ at $0\in \mathbb{C}^2$ linked to $\omega$ by an equisingular unfolding such that  the singularities reduction of $\omega'$ is exactly $E'$.
\end{Thm}

\noindent The last section is devoted to studying the relations between the moduli of a foliation and its invariant analytical curves, the so called separatrix. We use the existence theorem to prove the following

\begin{Thm}
Let $\FF_0$ be a non-dicritical generalized curve at $0\in \mathbb{C}^2$. For any curve $S\subset \mathbb{C}^2,0$ topologically equivalent to $\textup{Sep}(\FF_0)$, there exists a germ of foliation $\FF$ at $0\in\mathbb{C}^2$ topologically equivalent to $\FF$ with $\textup{Sep}(\FF)=S$.
\end{Thm}
\noindent The previous theorem can be expressed in terms of moduli spaces: the natural map $\mathtt{M}(\FF)\rightarrow \mathtt{M}(\text{Sep}(\FF))$ is onto, where $\mathtt{M}(.)$ refers to the moduli space.

\section{The categories $\text{Glu}_n(\MM,Z,\UU)$.}\label{part1}
The aim of this section is to describe a family of sub-categories in
{\it the marked trees category} (\ref{ArMr}).  These categories are built thanks to a {\it gluing process} (\ref{treeglu}) over a fixed tree and present some good properties of computation (\ref{Calc}). 

\subsection{The marked trees category.}\label{ArMr}

A blowing-up process over $\mathbb{C}^2$ is a commutative diagram 
\begin{equation}\label{ProEcla}
\begin{array}{ccccccccccccc}
\MM^h & \stackrel{E^h}{\rightarrow} & \ldots & \rightarrow &\MM^{j} & \stackrel{E^j}{\rightarrow} & \MM^{j-1} &
\rightarrow & \ldots & \stackrel{E^1}{\rightarrow} & \MM^{0} & = &
\mathbb{C}^2 \\
\bigcup &  &  &  &\bigcup &  &\bigcup  &
 &  &  & \bigcup & & \\
\Sigma^h & \rightarrow & \ldots & \rightarrow &\Sigma^{j} & \rightarrow & \Sigma^{j-1} &
\rightarrow & \ldots & \rightarrow & \Sigma^{0} & = & \{0\}\\
\bigcup &  &  &  &\bigcup &  &\bigcup  &
 &  &  & \bigcup & & \\
S^h & \rightarrow & \ldots & \rightarrow &S^{j} & \rightarrow & S^{j-1} &
\rightarrow & \ldots & \rightarrow & S^{0} & = & \{0\}\\
\end{array}\end{equation}
where $\MM^j$ is an analytical two-dimensional manifold; $\Sigma^j$ is a finite subset of the exceptional divsor $\DD^j:=(E^1\circ\cdots \circ E^j)^{-1}(S^0)$; $E^{j+1}$ is the standard blowing-up centered in $S_j$. The set $\Sigma_j$ is called the set of {\it singular points}. The set of irreducible components of $\DD^j$ is denoted by $\text{Comp}(\DD^j)$. The integer $h$ is called {\it height} of the blowing-up process and
$\left(\MM^h,\DD^h,\Sigma^h\right)$ {\it the top} of the process. The composed map $E_h:=E^1\circ\cdots \circ E^h$ is called the {\it total morphism} of the process.

\noindent Blowing-up process appears naturally in the desingularization theory. This article focusses on an isolated singularity of holomorphic foliation $\FF$ in $\mathbb{C}^2$. In this context, the Seidenberg's theorem \cite{reduction}\cite{MM} claims that there exists a blowing-up process $\left(E^j,\MM^j,\Sigma^{j},S^{j}\right)_{j=1\ldots h}$ reducing the singularities of $\FF$. A singularity of foliation is {\it reduced} when given by a holomorphic $1$-form with linear part
$$\lambda x\dd y +\beta y\dd x, \quad \beta\neq 0\quad\frac{\lambda}{\beta}\not\in\mathbb{Q}_{< 0}.$$
In the previous diagram, $\Sigma^j$ refers to the singularities of ${E^j}^*\FF$ and $S^j$ to the non-reduced ones. The foliation ${E^h}^*\FF$ requires  all its singularities to be reduced. 

\noindent More generally, we call {\it tree} a triplet
$(\MM,\DD,\Sigma)$, where $\MM$ is a two dimensional holomorphic manifold germ with a closed normal crossing hypersurface $\DD$. Each irreducible component of $\DD$ is biholomorphic to $\mathbb{P}^1$; $\Sigma$ is a finite subset of $\DD$ that contains the singular locus of $\DD$. Let us denote by $\left<\MM,\DD\right>$ the matrix
$$\left[<D,D'>\right]_{D,D'\in \text{Comp}(\DD)}$$ where $<D,D'>$ is the intersection number of the components $D$ and $D'$ . It's called {\it the tree intersection matrix}.

\bigskip
\noindent The {\it marked tree notion} was introduced in \cite{Seguy} in order to compare the foliation semi-local invariants. Let $(\MM^h_0,\DD^h_0,\Sigma^h_0)$ be the top of a blowing-up process.
We call {\it indexation} of $(\MM,\DD,\Sigma)$ related to $(\MM^h_0,\DD^h_0,\Sigma^h_0)$ a couple of two bijections $$\sigma:\Sigma^h_0\xrightarrow{\sim} \Sigma \quad \kappa:\text{Comp}(\DD^h_0)\xrightarrow{\sim}\text{Comp}(\DD)$$
such that for any irreducible component $D\subset \DD^h_0$, $\sigma(D\cap\Sigma^h_0)=\kappa(D)\cap\Sigma$.
Moreover, let us denote by $\Sigma\check{\times}\DD$ the following set
$$\left\{(s,D)\ \big|\ s\in D\right\}\subset \Sigma\times\text{Comp}(\DD).$$
We call  {\it complete closed paths system of $(\MM,\DD,\Sigma)$} a collection $$\Gamma=\fami{\gamma_{s,D}}_{(s,D)\in \Sigma\check{\times}\DD}$$ of closed paths in $\DD\backslash \Sigma$ such that, for any component $D$, the paths $\gamma_{s,D}$ generate the fundamental group of $D\backslash \Sigma$. The group $\text{Homeo}^+(\DD,\Sigma)$ of orientation preserving homeomorphisms, which fix each point of $\Sigma$, acts on the set of complete closed paths systems. Two complete closed paths systems $\Gamma$ and $\Lambda$ are isotopic when there exists a path $t\mapsto h_t$, $t\in [0,1]$ in $\text{Homeo}^+(\DD,\Sigma)$ such that $h_0=\text{Id}$ and $\Lambda=h_1^*\Gamma$. 

\begin{definition} A marking of $(\MM,\DD,\Sigma)$ related to $(\MM^h_0,\DD^h_0,\Sigma^h_0)$ is a triplet $(\sigma,\kappa,\Gamma)$ with an indexation $(\sigma,\kappa)$ and a complete path system isotopy class $\Gamma$ such that $\kappa$ conjugates the intersection matrix:
$$\left[<\kappa(D),\kappa(D')>\right]_{D,D'\in \text{Comp}(\DD^h_0)}=\left<\MM^h_0,\DD^h_0\right>.$$
\end{definition}
\noindent The {\it marked weighted dual tree}, denoted by $\mathbb{A}^*[\MM,\DD,\Sigma]$, is a couple of data: the weighted dual graph whose incidence matrix is the intersection matrix and the tree indexation. Two trees $(\MM,\DD,\Sigma)$ and
$(\mathcal{N},\mathcal{E},\Delta)$ with marking $(\sigma,\kappa,\Gamma)$ and $(\rho,\theta,\Lambda)$ are {\it conjuguated} when there is a biholormphism germ $H$ defined on a neighborhood of $\DD$ such that:
\begin{enumerate}
\item $H(\DD)=\mathcal{E},\qquad H(\Sigma)=\Delta;$
\item $H^*\sigma=\rho, \qquad H^*\kappa=\theta, \qquad H^*\Gamma=\Lambda.$
\end{enumerate}
The assumption $(2)$ is called {\it marked conjugation compatibility}. Let us denote by $\mathfrak{A}\left(\MM,\DD,\Sigma\right)$ the category of marked trees by  $\left(\MM,\DD,\Sigma\right)$ with marked compatible conjugation as arrows. It's not hard to prove the following proposition by induction on the height:

\begin{prop}\label{Cat} Let $(\NN,\mathcal{E},\Delta)\in \mathfrak{A}\left(\MM,\DD,\Sigma\right)$. Then there exists a blowing -up process whose marked top is conjugated to $(\NN,\mathcal{E},\Delta)$.
\end{prop}
\noindent The above proposition (\ref{Cat}) allows us to extend in an easy way all natural invariants of blowing-up process to the category of marked trees. Afterwards, we are going to define {\it the component multiplicity}. 

\noindent Let us denote by $\OO(\MM^h_0)$ the sheaf of germs of holomorphic function over $\MM^h_0$. Let $i_{\DD^h_0}$ be the divisor inclusion  $\DD^h_0\subset \MM^h_0$. We define
$$\OO_{\MM^h_0}:=i_{\DD^h_0}^{-1}(\OO(\MM^h_0)).$$
The sheaf $\OO_{\MM^h_0}$ is the restriction of $\OO$ over $\DD^h_0$. Let $D$ be an irreducible component of the divisor. We denote $$I_D\subset \OO_{\MM^h_0}$$ the ideal subsheaf of germs of function vanishing along $D$. Let us consider the subsheaf $\mathfrak{M}\subset \OO_{\MM^h_0}$ pull-back of the maximal ideal at $0\in\mathbb{C}^2$. The sheaf $\mathfrak{M}$ is locally free and generated by two global sections. Futhermore, we have the following decomposition 
\begin{equation*}
\mathfrak{M}=\prod_{D\in Comp(\DD)} I_{D}^{\nu(D)}.
\end{equation*}
Here, $\nu(D)$ is called {\it the multiplicity of the component $D$}. One can see that it is well determined by the intersection matrix. When  $\left(\MM,\DD,\Sigma\right)$ is a general tree marked by $(\MM^h_0,\DD^h_0,\Sigma^h_0)$, the multiplicity of $D$ is naturally the multiplicity of the component associated by the indexation in $\DD^h_0$.

\subsection{The sheaves $\GG^n_Z$, $n\geq 0$.}
\begin{quote}
{\bf From now on, we fix an element $(\MM,\DD,\Sigma)$ in $\mathfrak{A}\left(\MM_0,\DD^h_0,\Sigma^h_0\right)$}.\end{quote}
 In order to define some sub-categories of $\mathfrak{A}\left(\MM_0,\DD^h_0,\Sigma^h_0\right)$ we are interested in, we introduce a
family of sheaves of group over $\DD$. This construction leads to a key property (\ref{cle}). 

\noindent In order to get through a technical difficulty, which appears in the final induction (\ref{induction}), the tree is enhanced with a {\it cross}: let $\MM^h$ be the top of a blowing-up process with $\MM^h\stackrel{\phi}{\simeq}\MM$ given by the proposition (\ref{Cat}); let $E_h$ be the total morphism of the process and $E:=\phi\circ E_h$.
\begin{definition}[Cross]\label {Croix}
A cross on $\MM$ is the strict transform $Z=E^*Z_0$ of a single $Z_0=\fami{Z_1}$ or of a couple $Z_0=\fami{Z_1,Z_2}$ of germs of smooth transversal curves at the origin of $\mathbb{C}^2$.
\end{definition}
\noindent Throughout this article, we will often have to describe objects in coordinates. {\it Adapted coordinates} will always refer to $(x,y)$ local coordinates such that
\begin{itemize}
\item in the neighborhood of a regular point of $\DD$, $\{x=0\}$ is a local equation of $\DD$;
\item in the neighborhood of a singular point of $\DD$, $\{xy=0\}$ is a local equation of $\DD$;
\item in the neighborhood of the cross, $\{x=0\}$ is a local equation of $\DD$ and $\{y=0\}$ an equation of $Z$.
\end{itemize}

\noindent We consider
$\text{Aut}(\mathcal{M},Z)$ the groups sheaf over $\DD$ of germs of automorphism defined in a neighborhood of $\DD$ such that
$$ \Phi_{|\DD}=\textup{Id},\qquad \Phi_{|Z}=\textup{Id}.$$

\noindent Let us have a close look at the form of the  
$\textup{Aut}(\mathcal{M},Z)$ sections. At a regular point $c$ of $\DD\cup Z$
in an adapted coordinates system $(x,y)$ , the stack
$\textup{Aut}(\mathcal{M},Z)_c$ is the set of germs 
$(x,y)\mapsto \left(x(\alpha+A),y+xB\right)$
where $A$ and $B$ belong to $\mathbb{C}\{x,y\}$, $A(0,0)=0$ and $\alpha\in\mathbb{C}^*$. At a singular point $s$ of $\DD\cup Z$ in an adapted coordinates system, the stack $\textup{Aut}(\mathcal{M},Z)_s$ is the set of germs
$
(x,y)\mapsto \left(x(1+yA), y(1+xB)\right)
$ where $A$ and $B$ belong to $\mathbb{C}\{x,y\}$. 

\bigskip
\noindent We will be naturally led to look at the infinitesimal neighborhood of the divisor. To take care of the cross, we consider a filtration of  $\OO_\MM$ defined by
$\mathfrak{M}^n_Z:=I_Z\cdot\mathfrak{M}^n,\  n\geq 1$.
In the same way, we denote $\mathfrak{I}_Z\subset\OO_\MM$ the sheaf 
\begin{eqnarray*}
\mathfrak{I}_Z&:=&I_Z\cdot\prod_{D\in Comp(\DD)} I_{D}.
\end{eqnarray*}
\begin{definition}[Infinitesimal crossed tree] We call
$n^{th}$ infinitesimal crossed tree the analytical space
$$\MM^{[n],Z}:=\big( \DD, \OO_\MM\left/{\mathfrak{M}^n_Z}\right. \big).$$
The neighborhood of order $0$ is 
$\MM^{[0],Z}:=\big( \DD, \OO_\MM\left/{\mathfrak{I}_Z}\right. \big).$
We also consider the following ringed spaces:
\begin{equation*}
\MM^{\underline{n},Z} :=\big( \DD,\mathfrak{I}_Z\left/\mathfrak{I}_Z\mathfrak{M}^n_Z\right. \big) \quad
\MM^{\underline{0},Z} :=\big( \DD,\mathfrak{I}_Z\left/\mathfrak{I}^2_Z\right. \big). 
\end{equation*}
\end{definition} 
\noindent The sequence of canonical imbeddings 
$$\ldots \MM^{[p],Z} \hookrightarrow  \MM^{[p-1],Z}  \hookrightarrow \ldots  \hookrightarrow \MM^{[1],Z} \hookrightarrow \MM^{[0],Z} \hookrightarrow \MM $$
induces a natural filtration of the sheaf
$\text{Aut}(\MM,Z)$:
\begin{definition} We denote by $\textup{Aut}_n(\mathcal{M},Z)$ the subsheaf of $\textup{Aut}(\mathcal{M},Z)$ of germs that coincide with $\textup{Id}$ when restricted to 
the infinitesimal neighborhood of order $n$.
\end{definition}
\noindent The sequence of normal inclusions of groups  
$$ \ldots \triangleleft
\textup{Aut}_p(\mathcal{M},Z)\triangleleft\textup{Aut}_{p-1}(\mathcal{M},Z)\triangleleft\ldots\triangleleft\textup{Aut}_0(\mathcal{M},Z)
\triangleleft\textup{Aut}(\mathcal{M},Z)$$
points out the existence of a main part function for the filtration above. To describe it, we study the form of 
$\textup{Aut}_n(\mathcal{M},Z)$ sections. 

\noindent {\bf At a regular point $c$ of $\DD\cup Z$:} let $p$ be the multiplicity of the $c$ component. In an adapted coordinates system, the elements of
$\textup{Aut}_n(\mathcal{M},Z)_c$ can be written
$\phi(x,y)=\left(x+x^{pn}A,y+x^{pn}B\right)$,
where $A,B$ belong to
$\mathbb{C}\{x,y\}$. Let $\mathcal{J}_n$ be the function defined by
$$\phi\in\textup{Aut}_n(\mathcal{M},Z)_c\stackrel{\mathcal{J}_n}{\longmapsto} x^{pn-1}A(x,y)\in\left(\OO_{\MM^{[n],Z}}\right)_c. $$
One can see that $\mathcal{J}_n$ is a morphism of groups that doesn't depend on the adapted coordinates.

\noindent {\bf At a singular point $s$ of $\DD$}: let $p$ and
$q$ be the multiplicities of the local components. The elements of $\textup{Aut}_n(\mathcal{M},Z)_s$ are these of the form $\phi(x,y)=\left(x+x^{pn}y^{qn}A,y+x^{pn}y^{qn}B\right)$.
In the same way, we define an intrisic group morphism by
\begin{equation*}
\phi\in \textup{Aut}_n(\mathcal{M},Z)_s\stackrel{\mathcal{J}_n}{\longmapsto} x^{pn-1}y^{qn-1}\left(yA(x,y)+xB(x,y)\right) \in
\left(\OO_{\MM^{[n],Z}}\right)_s .
\end{equation*}

\noindent {\bf At an attach point $z$ of $Z$}: the elements of $\textup{Aut}_n(\mathcal{M},Z)_z$ are of the form
$\phi(x,y)=\left(x+x^{pn}yA,y+x^{pn}yB\right)$
and the morphism is defined by
\begin{equation*}
\phi \in \textup{Aut}_n(\mathcal{M},Z)_z\stackrel{\mathcal{J}_n}{\longmapsto} x^{pn-1}\left(yA(x,y)+xB(x,y)\right)\in \left(\OO_{\MM^{[n],Z}}\right)_z.
\end{equation*}\noindent Finally, we get a morphism of sheaves defined by its previous local description
$\textup{Aut}_n(\mathcal{M},Z)\stackrel{\mathcal{J}_n}{\longmapsto}\OO_{\MM^{[n],Z}}.$
We have likewise a morphism of sheaves $\mathcal{J}_0$:
\begin{equation}\label{superj0}
\textup{Aut}_0(\mathcal{M},Z)\stackrel{\mathcal{J}_0}{\longmapsto} \OO_{\MM^{\underline{0},Z}}.
\end{equation}
\begin{definition} We denote $\GG^n_Z$ the subsheaf of
  $\textup{Aut}_n(\mathcal{M},Z)$ kernel of the morphism $\mathcal{J}_n$.
\end{definition}
\noindent Lemma (\ref{criterion}) gives an intrisic criterion for a germ to be a section of $\GG^n_Z$:

\begin{lemme}\label{criterion} Let $\phi$ a germ of section of
  $\textup{Aut}(\mathcal{M},Z)$. The following properties are equivalent:
\begin{enumerate}
\item $\phi$ is a section of $\GG^n_Z$.
\item $\phi$ is the identity restricted to $\MM^{[n],Z}$ and to $\MM^{\underline{n},Z}$.
\end{enumerate}
\end{lemme}
\noindent The following property shows the interest of the sheaves $\GG^n_Z$ we introduce: basically, these are sheaves of Lie groups associated to some sheaves of Lie algebras which are natural in our context.
\begin{propr}[Key property]\label{cle}
Let $X$ be a germ of vector field tangent to $\DD$ and to $Z$. Let $f$ be a germ of $\mathfrak{M}^n_Z$ section. Then the flow of $f\cdot X$ is a germ of $\GG^n_Z$ section.
\end{propr}
\begin{demo} The property can be read on the form of flow in coordinates. Let $(x,y)$ be local adapted coordinates. If $Y$ is a germ of vector field, we denote by $Y^{(k)}$ the  $k^{th}$ power of $Y$ as differential operator on $\OO_{0}^2$; for $x$ and $y$ small enough and $t\in[0,1]$, the flow of $Y$ at time $t$ can be expanded as
$$e^{(t)Y}=\sum_{k=0}^{\infty}\frac{t^k}{k!}Y^{(k)}(x,y).$$
In the case of a divisor singular point, an induction on $k$ shows that, there exists $A_k,B_k$, sections of $
\mathfrak{M}^{n+k}_Z$, such that
$Y^{(k)}(x,y)=(xA_k(x,y),yB_k(x,y))$.
Hence, at time $1$, the flow has the form
$(x,y)\longmapsto (x,y) +x^{pn}y^{qn}(xA,yB)$ where $p$ and $q$ are the multiplicities of the local components of $\DD$. The latter automorphism is a $\GG^n_Z$ section.
\end{demo}

\subsection{The tree gluing.}\label{treeglu}

Thanks to the sheaf $\text{Aut}(\MM,Z)$, we are going to introduce a process called
{\it gluing} on $\MM$. This construction will allow us to define a large class of trees with same divisor analytical type. These trees will inherit a canonical marking and a cross.

\subsubsection{Distinguished covering.}\label{distingue}
Let us define a particular type of open covering of the divisor. Open sets of that covering will play the role of gluing "bricks "

\medskip
\noindent Let $\UU=\left\{U_i\right\}_{i \in \mathbb{I}=\mathbb{I}_0\cup \mathbb{I}_1}$ be covering of $\DD$ constitued of two kinds of open sets:
\begin{itemize}
\item if $i$ belongs to $\mathbb{I}_0$, $U_i$ is the trace on $\DD$ of a neighborhood of a unique singular.
\item if $i$ belongs to $\mathbb{I}_1$, $U_i$ is an irreducible component of $\DD$ deprived of the singular points of $\DD$.
\end{itemize}

\begin{definition}
Every such covering is called distinguished when there is no $3$-intersection.
\end{definition}
\noindent Distinguished coverings contain Stein open sets having fundamental systems of Stein neighborhood. From now on, a covering
denoted by $\UU$ will always supposed to be distinguished. The spaces $\ZZ^0\left(\UU,\GG\right)$ and $\ZZ^1\left(\UU,\GG\right)$
are the sets of $0$-cocycles and $1$-cocycles in the sense of Cech
for the sheaf $\GG$ and the covering
$\UU$. Let us consider
$$\mathbb{I}_0\check{\times}\mathbb{I}_1=\left\{(i,j)\in \mathbb{I}_0\times\mathbb{I}_1|U_i\cap U_j \neq \emptyset \right\}.$$
We define 
$\tilde{\ZZ}^1:=\displaystyle\prod_{(i,j)\in\mathbb{I}_0\check{\times}\mathbb{I}_1}\ZZ^0\left(U_i\cap U_j,\textup{Aut}(\mathcal{M},Z)\right)$.
Since distinguished covering doesn't have any $3$-intersection,  $\tilde{\ZZ}^1$ and $\ZZ^1 \big(\UU,\textup{Aut}(\mathcal{M},Z)\big)$ are isomorphic. Hence, if no confusion is possible, we will keep on denoting $\ZZ^1$ the space $\tilde{\ZZ}^1$.

\subsubsection{Gluing.}\label{collage}
Thanks to distinguished covering, we are able to glue the open sets of that covering by identifying points according to a $1$-cocycle in $\textup{Aut}(\mathcal{M},Z)$. Let 
$\coup{\phi_{ij}}$ be a $1$-cocycle in $\ZZ^1\left(\UU,\textup{Aut}(\mathcal{M},Z)\right)$. We define:
$$\MM[\phi_{ij}]= \bigcup_i \mathcal{U}_i\times\{i\}
/_{\big(\{x\}\times\{i\}\sim \{\phi_{ij}(x)\}\times\{j\}\big)}, $$
 where $\UU_i$ is a neigborhood of $U_i$ in $\MM$ such that $\phi_{ij}$ becomes an automorphism of $\MM$ along $U_i\cap U_j$. The latter automorphism is well defined since it coincides with $\text{Id}$ along the divisor. The manifold we get comes with an embedding 
\begin{equation}\label{inj}
\DD \hookrightarrow \MM[\phi_{ij}]
\end{equation} 
whose image is denoted by $\DD[\phi_{ij}]$ and $\MM[\phi_{ij}]$ is considered as a germ of neighborhood of $\DD[\phi_{ij}]$.

\begin{definition}
The manifold germ $\MM[\phi_{ij}]$ is called gluing of $\MM$ along $\UU$ by the cocycle $\coup{\phi_{ij}}$. 
\end{definition}
\noindent The gluing of a marked crossed tree comes naturally with a marking and a cross: the marking is the direct image by the embedding (\ref{inj}) of the marking $\sigma$; the cross is the direct image of $Z$ by the quotient map for gluing relation. Such a tree, marking and cross are respectively denoted by
\begin{eqnarray*}
\left(\MM[\phi_{ij}],\DD[\phi_{ij}],\Sigma[\phi_{ij}]\right),\ \sigma[\phi_{ij}],\text{ and } Z[\phi_{ij}].
\end{eqnarray*}
In the same way, the direct image of the covering $\UU$ by the quotient map is a distinguished covering of the new tree and is denoted by $\UU[\phi_{ij}]$. 

\medskip
\noindent We associate to any gluing the data of morphisms on infinitesimal neigborhood generalizing the embedding (\ref{inj}). The intrisic description (\ref{criterion}) of $\GG_n^Z$ sections reveals the following property

\begin{propr}\label{superinj} Let $n$ be an integer and $\NN=\MM[\phi_{ij}]$ be a gluing of $\MM$ by a cocycle in $\ZZ^1(\UU,\GG^n_Z)$. Then the canonical analytical embeddings 
\begin{eqnarray*}
\rho^{[n]}_{\NN}:&\MM^{[n],Z}\hookrightarrow \NN, \\
\rho^{\underline{n}}_{\NN}:&\MM^{\underline{n},Z}\hookrightarrow \NN
\end{eqnarray*}
have for respective images $\NN^{[n],Z[\phi_{ij}]}$ and $\NN^{\underline{n},Z[\phi_{ij}]}$.
\end{propr}

\subsection{The $\text{Glu}_n(\MM,Z,\UU)$ categories.}
Let $p$ be an integer. Let us consider the marked crossed tree built by a succession of gluings
\begin{equation}\label{var1}
\MM[\phi^1_{ij}][\phi^2_{ij}][\ldots][\phi^p_{ij}]
\end{equation}
where
\begin{itemize}
\item $\coup{\phi^1_{ij}}$ is a $1$-cocyle of $\GG^n_Z$;
\item for $k=2,\ldots,p$, $\coup{\phi^k_{ij}}\in\ZZ^1\left(\UU[\phi^1_{ij}]\ldots[\phi^{k-1}_{ij}],\GG^n_{Z[\phi^1_{ij}]\ldots[\phi^{k-1}_{ij}]}\right)$ with $$\GG^n_{Z[\phi^1_{ij}]\ldots[\phi^{k-1}_{ij}]}\subset \text{Aut}\left(\MM[\phi^1_{ij}]\ldots[\phi^{k-1}_{ij}],Z[\phi^1_{ij}]\ldots[\phi^{k-1}_{ij}]\right).$$ 
\end{itemize}
Following the (\ref{superinj}) lemma, we have canonical embeddings
\begin{eqnarray}
\rho^{[n]}_{\MM[\phi^1_{ij}][\phi^2_{ij}][\ldots][\phi^p_{ij}]}:&\MM^{[n],Z}\hookrightarrow \MM[\phi^1_{ij}][\phi^2_{ij}][\ldots][\phi^p_{ij}],\label{inj1} \\
\rho^{\underline{n}}_{\MM[\phi^1_{ij}][\phi^2_{ij}][\ldots][\phi^p_{ij}]}:&\MM^{\underline{n},Z}\hookrightarrow \MM[\phi^1_{ij}][\phi^2_{ij}][\ldots][\phi^p_{ij}]. \label{inj2}
\end{eqnarray}

\begin{definition}
The $\textup{Glu}_n(\MM,Z,\UU)$ category is the category whose objects are marked crossed trees built as (\ref{var1}) with the data of embeddings  (\ref{inj1}) and (\ref{inj2}). Arrows are biholomorphic germs that respect marking and commute with embeddings. If $\MM$ and $\NN$ are isomorphic in $\textup{Glu}_n(\MM,Z,\UU)$, we denote
$$\MM\stackrel{\mathtt{G}_n}{\simeq} \NN.$$ 
\end{definition}

\subsection{Computations in $\textup{Glu}_0(\MM,Z,\UU)$.}\label{Calc}
To compute in $\textup{Glu}_0(\MM,Z,\UU)$, we state three properties to manipulate gluings towards their defining cocycles. In view of intrisic description of $\GG_n^Z$ sections, it's easy to check the following property:
\begin{propr}\label{lell}Let $\NN$ and $\PP$ be in $\text{Glu}_n(\MM,Z,\UU)$
  defined respectively by the $1$-cocyles $\coup{\rho_{ij}}$ and $\coup{\gamma_{ij}}$.
The following properties are equivalent:
\begin{enumerate}
\item $\NN\stackrel{\mathtt{G}_n}{\simeq} \PP.$ 
\item There exists a $0$-cocyle
$\coup{\phi_i}$ in $\GG^n_Z$ such that
$\rho_{ij}=\phi_j \gamma_{ij}\phi_i^{-1}$.
\end{enumerate}
\end{propr}
\noindent The family of maps
$I_i:U_i\subset\MM \rightarrow \MM[\phi_{ij}]$
induces the following canonical isomorphism 
$$\zeta^0:\left\{
\begin{array}{ccc}
\ZZ^0 \big(\mathcal{U},\textup{Aut}(\mathcal{M},Z)\big) & \longrightarrow &
\ZZ^0\big(\mathcal{U}[\phi_{ij}],\textup{Aut}(\mathcal{M}[\phi_{ij}],Z[\phi_{ij}])\big)\\
\coup{\phi_i} & \longrightarrow & \coup{I_i\phi_i I_i^{-1}} \\
\end{array}\right. .$$
We are able to define such an isomorphism for $1$-cocycles thanks to the $\tilde{\ZZ}^1$ spaces (\ref{distingue}): we define $\zeta^1$ as the following isomorphism,
$$\begin{array}{ccc}
\tilde{\ZZ}^1 & \rightarrow &
\tilde{\ZZ}^1[\phi_{ij}]:=\prod_{(i,j)\in\mathbb{I}_0\check{\times}\mathbb{I}_1}\ZZ^0\left(U_i\cap U_j[\phi_{ij}],\textup{Aut}(\mathcal{M} [\phi_{ij}],Z[\phi_{ij}])\right) \\
\coup{\phi_{ij}} & \rightarrow & \coup{I_j\phi_{ij} I_j^{-1}} 
\end{array} $$
is an isomorphism and the checked map $\zeta^1$. Now, one can state the following useful property:

\begin{propr}\label{lell1}
Let $\coup{\phi_{ij}}$ and $\coup{\psi_{ij}}$ be in
$\ZZ^1\big(\UU,\textup{Aut}(\mathcal{M},Z)\big)$. Then
$$ \MM[\phi_{ij}][\zeta^1\psi_{ij}]\stackrel{\mathtt{G}}\simeq\MM[\phi_{ij}\psi_{ij}].$$
\end{propr}

\subsubsection{Stability property.}
In this subsection, we identify in the gluing category a class of isomorphic trees: basically, if the gluing cocycle is tangent enough, the glued tree is isomorphic to the initial one.

\begin{propr}[Stability]
\
\noindent For $n$ big enough, the image of
$$\textup{Glu}_n(\MM,Z,\UU)\longrightarrow \textup{Glu}_0(\MM,Z,\UU)$$ are trees isomorphic to $\MM$. 
\end{propr} 
\noindent This property can be stated in the following way: the map
$$H^1\left(\DD,\GG^n_Z\right)\longrightarrow H^1\left(\DD,\GG^0_Z\right)$$
is constant equal to $\left[Id\right]_{\GG^0_Z}$. First, let us establish the equivalent statement for the intial sheaf $\text{Aut}_n(\DD,Z).$
\begin{lemme}\label{triv}
There is an integer $\delta(n)$ bigger than $n$ such that
the image of the map
$$H^1(\DD,\textup{Aut}_{\delta(n)}(\MM,Z))\rightarrow H^1(\DD,\textup{Aut}_{n}(\MM,Z))$$
is trivial.
\end{lemme}
\begin{demo} The proof of a similary result can be found in \cite{MatSal}. Here, we reproduce the main arguments. One can supposed $\MM$ to be the top of a blowing-up process. Let $p\geq n$ and $\coup{\phi_{ij}}$ be an element of 
$Z^1\left(D,\GG^p_Z\right)$. There exists a germ of biholomorphism 
$\theta$ between $\MM[\phi_{ij}]$ and the top of blowing-up process $\MM'$. Let us denote by $Z'$ the induced cross on $\MM'$. The map $\theta$ induces an isomorphism $\theta^{[p],Z}$ between the infinitesimal crossed neighborhoods $\MM[\phi_{ij}]^{[p],Z[\phi_{ij}]}$ and $\MM'^{[p],Z'}$. 
Hence, $\theta^{[p],Z}\circ\rho^{[p],Z}$ identifies infinitesimal neighborhoods of top's blowing-up process. One can show that for $p=\delta(n)$ big enough, $\theta^{[p],Z}\circ\rho^{[p],Z}$ can be extended in $T$ as a biholomorphism of trees such that
$$ T^{[n],Z}=\theta^{[n],Z}\circ\rho^{[n],Z}.$$
Hence, $H=T^{-1}\circ\theta$ is a germ of biholomorphism between $\MM[\phi_{ij}]$ and
$\MM$ with $H^{[n],Z}\circ \rho^{[n],Z}=\text{Id}^{[n],Z}$. This is equivalent to the triviality of $\coup{\phi_{ij}}$ in $H^1(\DD,\text{Aut}_n(\MM,Z))$.
\end{demo}
\begin{coro}\label{stabcoro}
There is an integer $n_0$ such that for $n\geq n_0$ the image of
$$H^1(\DD,\GG^n_Z)\rightarrow H^1(\DD,\GG^0_Z)$$
is trivial.
\end{coro}
\begin{demo} Let $\coup{\phi_{ij}}$ be a $1$-cocycle in $\GG^{\delta(n)}_Z$. Its trivialisation in $\textup{Aut}_{n}(D,Z)$ can be written
\begin{equation}\label{coco}
\phi_{ij}=\phi_{i}\circ\phi_{j}^{-1}.
\end{equation}
Let us denote by $p_1$ and $p_2$ the attached points of the cross.  
The cross is the strict pulled-up of a couple of transversal curves at the origin in $\mathbb{C}^2$ which define a local system of coordinates. Near $p_1$ and $p_2$, there are two local adapted coordinates systems $(x_1,y_1)$ and $(x_2,y_2)$ such that $E$ takes the following forms:
$E(x_1,y_1)=(x_1,y_1x_1^{N_1})$ and
$E(x_2,y_2)=(x_2y_2^{N_2},y_2)$.
The components of the cocycle $\phi_1$ and $\phi_2$ defined near $p_1$ and $p_2$ can be expanded as:
\begin{eqnarray*}
\phi_1(x_1,y_1)=\left(x_1+x_1^ny_1U_1(x_1,y_1),y_1+x_1^ny_1V_1(x_1,y_1) \right) \\
\phi_2(x_2,y_2)=\left(x_2+y_2^nx_2U_2(x_2,y_2),y_2+y_2^nx_2V_2(x_2,y_2) \right).
\end{eqnarray*}
Let $\phi_0$ be the germ of biholomorphism near $0$ in $\mathbb{C}^2$ defined by
$$\phi_0(x,y)=(xe^{y^nU_2(0,y)},ye^{x^nV_1(x,0)}).$$ 
For $n$ big enough, $\phi_0$ can be raised in an automorphism $\phi$ that fixes each point of the divisor and the cross. Moreover, for any point $c$ different from $p_1$ and $p_2$, the evaluation through $\mathfrak{J}_0$ (see (\ref{superj0})) provides the equality
$\mathfrak{J}_0(\phi^{-1})_c\equiv 0$.
Now, the choice of $\phi$ ensures that 
$\mathfrak{J}_0(\phi_1)_{p_1} \equiv \mathfrak{J}_0(\phi)_{p_1}$ and that
$\mathfrak{J}_0(\phi_2)_{p_2}  \equiv \mathfrak{J}_0(\phi)_{p_2}$.
Finally, $\mathfrak{J}_0(\phi_i\circ\phi^{-1})\equiv 0$. Hence, the $0$-cocycle $\coup{\phi_i\circ \phi^{-1}}$ is a trivialisation of $\coup{\phi_{ij}}$ in $\GG^0_Z$.
\end{demo}
As mentioned before, the previous corollary is equivalent to the stability property.
\section{Cobordism in $\text{Glu}_n(\MM,Z,\UU)$.}\label{mod}

From now on, we assume the marked tree $\MM$ to be foliated by $\FF$. We are going to define a cobordism notion in order to detect on any
element of $\text{Glu}_0(\MM,Z,\UU)$ the existence of foliation linked to 
$\FF$ in the sense of equisingular unfolding. For a precise definition of equisingular unfolding, we refer to \cite{univ}.
We assume the singularities of $\FF$ to be reduced.
\begin{definition}[Cross adapted to $\FF$]
\noindent Let $Z$ be a cross on $\MM$. $Z$ is said to be adapted to 
$\FF$ when each $Z_i$ verifies at least one of the following properties:
\begin{enumerate}
\item $Z_i$ is a separatrix of $\FF$.
\item $Z_i$ is attached at a regular point of $\FF$. 
\end{enumerate}
In the latter case, $Z_i$ will be transversal to the foliation.
\end{definition}
\noindent Now, we consider two natural sheaves:
\begin{itemize}
\item the sheaf $\mathfrak{X}_{S,Z}$ is a sheaf over $\DD$. In each point of $\DD$, its fiber is the space of holomorphic vector field germs of $\MM$ that are tangent to 
$\DD$, to the separatrix and to the cross;
\item the sheaf $\mathfrak{X}_{\FF,Z}\subset\mathfrak{X}_{S,Z}$ is the subsheaf of germs tangent to the foliation.
\end{itemize}
From now on,  $e^{tX}$ refers to the flow of the vector field $X$ at time $t$. One notices that, if $X$ is a section of $\mathfrak{I}_Z\mathfrak{X}_{S,Z}$ then, for all $t\in \mathbb{C}$, $e^{tX}$ exists as germ and defines a section of $\text{Aut}(\MM,Z)$.
\begin{definition}[Elementary cobordism]\label{def1}
Let $\NN\in\textup{Glu}_0(\MM,Z,\UU)$. 

\noindent $\NN$ is elementary 
$\FF$-cobordant to $\MM$ when there exists a $1$-cocycle $\coup{T_{ij}}$
of
$Z^1\big(\UU,\mathfrak{I}_Z\mathfrak{X}_{\FF,Z}\big)$ such that
$\NN$ and $\MM[e^{T_{ij}}]$ are isomorphic in $\textup{Glu}_0(\MM,Z,\UU)$.
\end{definition}
\noindent Here, $\NN$ inherit a canonical foliation defined by $$\FF[e^{T_{ij}}]:=\coprod_i \FF|_{\mathcal{U}_i}\Big
/_{x\sim e^{T_{ij}}x}. $$
Moreover, a rigidity result of Grauert \cite{univ} implies that the deformation $t \rightarrow \MM_t=\MM[e^{(t)T_{ij}}],\ t\in\dbar$ carries an equisingular unfolding between $\FF$ and $\FF[e^{T_{ij}}]$ in the sense of \cite{univ}.

\begin{definition}[General cobordism]
Let $\NN$ be in $\textup{Glu}_0(\MM,Z,\UU)$. 
$\NN$ is said to be
$\FF$-cobordant to $\MM$ if there exists a finite sequence of $1$-cocyles
$\coup{T^k_{ij}}_{k=1,\ldots,N}$ such that the two following conditions are verified:
\begin{enumerate}
\item for any $p=0,\ldots,N-1$, let $\mathfrak{X}_{\FF_p,Z_p}$ be the sheaf over $\DD[e^{T^1_{ij}}][\cdots][e^{T^{p}_{ij}}]$ of germs of vector field tangent to the following foliation and cross $$\FF_p=\FF[e^{T^1_{ij}}][\cdots][e^{T^{p}_{ij}}],\quad Z_p=Z[e^{T^1_{ij}}][\cdots][e^{T^{p}_{ij}}].$$ We assume $\coup{T^{p+1}_{ij}}$ is a $1$-cocycle with values in  $\mathfrak{X}_{\FF_p,Z_p}$.

\item $\NN\stackrel{\mathtt{G}_0}{\simeq}\MM[e^{T^1_{ij}}][\cdots][e^{T^N_{ij}}]$.
\end{enumerate}
We summarize this definition with the following notation:
$$
 \MM\fleche{\FF_1,Z_1}\MM_2\fleche{ \FF_2,Z_2}\cdots \fleche{ \FF_{N-1},Z_{N-1}}\MM_N\stackrel{\mathtt{G}_0}{\simeq} \NN
$$
\end{definition}

\subsection{Conctruction of cobordism.}
A foliation $\FF$ is a  non-dicritical generalized curve 
when each reduced singularity has two non-vanishing eigenvalues and when the divisor is invariant.

\begin{prop}\label{cobordisme0}
Let $\FF$ be a non-dicritical generalized curve on a marked tree $\MM$ crossed by $Z$. Any element in $\textup{Glu}_0(\MM,Z,\UU)$ is $\FF$-cobordant to $\MM$.
\end{prop}
\noindent The proof points out three different steps : first, we establish the result at an infinitesimal level, then for the sub-category $\textup{Glu}_1(\MM,Z,\UU)$ and finally, thanks to an induction on the height of the tree, for the  $\textup{Glu}_0(\MM,\ZZ,\UU)$ category.

\subsubsection{Step 1: the infinitesimal level.}
Let $H$ be a biholomorphism between $\MM$ and the top of a blowing-up process $\MM^h$. The induced foliation $\FF^h$ on $\MM^h$ is defined by a germ of holomorphic $1$-form $\omega$ at the origin of $\mathbb{C}^2$. We define $E:=E_h\circ H$ where $E_h$ is the total morphism of the process (\ref{ArMr}). The global
$1$-form $E^{*}\omega$ defines 
a morphism of sheaves
$$X\in \mathfrak{X}_{S,Z}\rightarrow E^*(\omega)(X) \in \OO_\MM.$$
Let a reduced equation of the separatrix be denoted by $f$.
\begin{lemme}\label{suisui2} 
There exists an exact sequence of sheaves
$$0 \longrightarrow \mathfrak{M}^n_Z\mathfrak{X}_{\FF,Z} \longrightarrow
\mathfrak{M}^n_Z\mathfrak{X}_{S,Z}\xrightarrow{E^*\omega(.)} \mathfrak{M}^n_Z
\left(f\circ E\right)\longrightarrow 0. $$
Here, $\left(f\circ E\right)$ is the sheaf of ideals generated by $f\circ E$ in $\OO_\MM$.
\end{lemme}
\begin{demo} Let us show that there exists an exact sequence of sheaves
\begin{equation}\label{suisui}
0 \longrightarrow \mathfrak{X}_{\FF,Z} \longrightarrow
\mathfrak{X}_{S,Z}\xrightarrow{E^*\omega(.)}\left(f\circ
  E\right)\longrightarrow 0.
\end{equation}
The exactness of the first part in (\ref{suisui}) is obvious since $\XX_{\FF,Z}=\text{ker}\left(E^*\omega\right)$. Let us compute the image of $E^*\omega(.)$. As $\FF$ is a non-dicritical generalized curve, near a divisor singular point $s$ there exists local coordinates such that $f\circ E =Ux^{p+1}y^{q+1}$ and
$E^*\omega=Vx^py^q(\lambda x(1+A)\dd y + y(1+B)\dd x)$, where $U$ and $V$
  are unities. Let $g$ be any element of $(\OO_\MM)_s$. The vector field
$X=yg\frac{U}{V}\partial y$ belongs to $(\mathfrak{X}_{S,Z})_s$ and
$E^*\omega(X)=gf\circ E$. At a regular point $c$ of the foliation, there exists coordinates such that
$f\circ E =Ux^{p+1}$
and $E^*\omega=Vx^p\dd x$, where $U$ and $V$
 are unities. The vector field 
$X=xg\frac{U}{V}\partial x$ belongs to $(\mathfrak{X}_{S,Z})_c$ and verifies $E^*\omega(X)=gf\circ E$ for any germ of function $g$. The other cases can be studied in the same way. Now, the sheaf $\mathfrak{M}^n_Z$ is locally principal. Hence, the sequence (\ref{suisui}) multiplied by $\mathfrak{M}^n_Z$ remains exact.
\end{demo}
\begin{lemme}\label{lemlem} For any $n\geq 1$,
$
H^1(\DD,\mathfrak{M}^n_Z)=0.$
\end{lemme}
\begin{demo}  The long exact sequence associated to the short exact one
$0 \rightarrow \mathfrak{M}^n_Z \rightarrow \mathfrak{M}^n \rightarrow\mathfrak{M}^n/\mathfrak{M}^n_Z \rightarrow 0$ is written (\cite{Godement})
$$\cdots \rightarrow H^0(\DD,\mathfrak{M}^n) \stackrel{\delta}{\rightarrow}H^0(\DD,\mathfrak{M}^n/\mathfrak{M}^n_Z)\rightarrow H^1(\DD,\mathfrak{M}^n_Z)\rightarrow H^1(\DD,\mathfrak{M}^n)\cdots.$$
The sheaf $\mathfrak{M}^n$ is generated by its global sections. Hence $
H^1(\DD,\mathfrak{M}^n)=0$ since $H^1(\DD,\OO_\MM)=0$ \cite{univ}. In order to conclude, it remains to prove that $\delta$ is onto. Out of the attached points $p_1$ and $p_2$ of $Z$, the sheaf $I_{Z}$ coincides with
$\OO_\MM$, therefore the fiber of $\mathfrak{M}^n/\mathfrak{M}^n_Z$ is trivial. Let us consider the 
local coordinates systems $(x_1,y_1)$ and $(x_2,y_2)$ near $p_1$ and $p_2$ introduced in the proof of (\ref{stabcoro}). We have
$$\left(\mathfrak{M}^n/\mathfrak{M}^n_Z\right)_{p_1}\simeq x_1^n\mathbb{C}\{x_1,y_1\}/x_1^ny_1\mathbb{C}\{x_1,y_1\}\simeq x_1^n\mathbb{C}\{x_1\}.$$
In the same way, $\left(\mathfrak{M}^n/\mathfrak{M}^n_Z\right)_{p_2}$ is isomorphic to  $y_2^n\mathbb{C}\{y_2\} $.
Hence, the space of $\mathfrak{M}^n/\mathfrak{M}^n_Z$ global sections is identified to
$ x_1^n\mathbb{C}\{x_1\}\bigoplus y_2^n\mathbb{C}\{y_2\}$. Let $S=x_1^na_1(x_1)\oplus y_2^na_2(y_2)$ be in the previous set. The germ of function
defined by $s(x,y)=x^na_1(x)+y^na_2(y)$ induces a global section $s\circ E$ of $\mathfrak{M}^n$. Now, the map $E$ takes the following form in coordinates:
$E(x_1,y_1)=(x_1,y_1x_1^{N_1})$ and $E(x_2,y_2)=(x_2y_2^{N_2},y_2)$. Hence $s\circ E$ verifies the following equalities 
\begin{eqnarray*}
(s\circ E)_{p_1}&\equiv & x_1^na_1(x_1)+y_1^nx_1^{nN_1}a_2(y_1x_1^{n_1})\equiv x_1^na_1(x_1)\in\left(\mathfrak{M}^n/\mathfrak{M}^n_Z\right)_{p_1}\\
(s\circ E)_{p_2}&\equiv & x_2^ny_2^{nN_2}a_1(x_2y_2^{N_2})+y_2^na_2(y_2)\equiv y_2^na_2(x_2)\in\left(\mathfrak{M}^n/\mathfrak{M}^n_Z\right)_{p_2},
\end{eqnarray*}
that means $\delta(s\circ E)=S$ and ends the proof.
\end{demo}

\noindent Now, we can state the infinitesimal equivalent of the property (\ref{cobordisme0}):
\begin{prop}[Infinitesimal cobordism]\label{lem23}\ 

\noindent The map
$$H^1(\DD,\mathfrak{M}^n_Z\mathfrak{X}_{\FF,Z})\longrightarrow
H^1(\DD,\mathfrak{M}^n_Z\mathfrak{X}_{S,Z})$$ is onto.
\end{prop}
\begin{demo} The long exact sequence associated to the short one in the lemma (\ref{suisui2}) is 
\begin{equation}\label{suitexa}
\cdots\rightarrow H^1(\DD,\mathfrak{M}^n_Z\mathfrak{X}_{\FF,Z})\rightarrow
H^1(\DD,\mathfrak{M}^n_Z\mathfrak{X}_{S,Z})\rightarrow
H^1(\DD,\mathfrak{M}^n_Z(f\circ E))\rightarrow\cdots
\end{equation}
Now, the lemma (\ref{lemlem}) ensures that $H^1(\mathfrak{M}^n_Z(f\circ E))$ is trivial,
hence the last term of (\ref{suitexa}) vanishes. \end{demo}

\subsubsection{Step 2: cobordism in $\textup{Glu}_1(\MM,Z,\UU)$.}
This section is devoted to prove the following proposition:
\begin{prop}\label{Glulg1}
Let $\FF$ be a non-dicritical generalized curve on a marked tree $\MM$ crossed by $Z$.
Any element of $\textup{Glu}_1(\MM,Z,\UU)$ is elementary $\FF$-cobordant to $\MM$.
\end{prop}
\noindent The proof consists in getting the cobordism on a infinitesimal neighborhood of order big enough to apply the stability property. 
\begin{lemme}[Infinitesimal cobordism of order $n$]
\

\noindent Let $\coup{\phi_{ij}}$ in
$\ZZ^1(\UU,\GG_Z^1)$. For any $n\geq 1$,  there exists $\coup{T_{ij}}$
in $\ZZ^1\big(\UU,\mathfrak{M}_Z\mathfrak{X}_{\FF,Z}\big)$, $\coup{\phi_i}$
in $\ZZ^0\left(\UU,\GG_Z^1\right)$ and $\coup{\tilde{\phi}_{ij}}$
in $\ZZ^1(\UU,\GG_Z^n)$ such that

\begin{equation}\label{rel1}
\phi_j^{-1}\circ\phi_{ij}\circ\phi_i=e^{T_{ij}}\circ\tilde{\phi}_{ij}.
\end{equation}

\end{lemme}
\begin{demo} It is an induction on the integer $n$. Let us assume the property to be true at rank $n$. Let $\coup{\tilde{\phi}_{ij}}\in\ZZ^1(\UU,\GG_Z^n)$ be given by the induction hypothesis. In local adapted coordinates, one can write
$\tilde{\phi_{ij}}(x_{ij},y_{ij})=\textup{Id}+(A_{ij},B_{ij})$.
Let $X_{ij}$ be the germ of vector field defined by $X_{ij}=A_{ij}\partial_{ x_{ij}} + B_{ij}\partial_{ y_{ij}}$.
The family
$\coup{\tilde{X}_{ij}}$ is a $1$-cocycle with values in $\mathfrak{M}^n_Z\mathfrak{X}_{S,Z}$. The lemma (\ref{lem23}) ensures the existence of a $0$-cocycle $\coup{\tilde{X}_i}$ in $Z^0\left(\UU,\mathfrak{M}^n_Z\mathfrak{X}_{S,Z}\right)$  and a $1$-cocycle $\coup{\tilde{T}_{ij}}$ in $Z^1\left(\UU,\mathfrak{M}^n_Z\mathfrak{X}_{\FF,Z}\right)$ such that
$
\tilde{X}_{ij}=\tilde{X}_j-\tilde{X}_i+\tilde{T}_{ij}
$.
By expanding the flow $e^{\tilde{X}_{ij}}$, one can see that
$$\phi_{ij}^1:=e^{-\tilde{X}_{ij}}\circ\tilde{\phi}_{ij}=(\text{Id}-\tilde{X}_{ij}+\cdots)\circ(\text{Id}+\tilde{X}_{ij}+\cdots)\in \GG^{n+1}_Z.$$
In view of the previous relations, we find
\begin{eqnarray*}
e^{-\tilde{X}_j}\circ\phi_j^{-1}\circ\phi_{{ij}}\circ \phi_i\circ e^{\tilde{X}_i}&=&e^{-\tilde{X}_j}\circ e^{T_{ij}}\circ\tilde{\phi}_{ij}\circ e^{\tilde{X}_i}\\
&=&e^{-\tilde{X}_j}\circ e^{T_{ij}}\circ e^{\tilde{X}_{ij}}\circ e^{\tilde{X_i}}\circ \phi_{ij}^1\circ\left[\phi_{ij}^1,e^{\tilde{X}_i}\right]
\end{eqnarray*}
where $[a,b]=a^{-1}b^{-1}ab$.
Now, one can find in \cite{spivak} the following result:
\begin{souslemme}[Campbell-Hausdorff formula]
\ 

\noindent Let $X$ and $Y$ be two germs of vector field vanishing along the divisor. Then
there exists a formal vector field  
$\rho(X,Y)=\sum_{k=1}^{\infty}\rho_k\left(X,Y\right)$
such that
$$e^{\rho(X,Y)}=e^{X}\circ e^Y.$$
The sum is convergent for the Krull topology in the space of formal series. Moreover, the first terms of the series are
\begin{equation}\label{ro2}
\rho\left(X,Y\right)=X+Y+\fraction{1}{2}\left[X,Y\right]+\cdots.
\end{equation}
\end{souslemme}
\noindent Thanks to the above lemma, we get the following relation \begin{equation*}
e^{-\tilde{X}_j}\circ e^{T_{ij}}\circ e^{\tilde{X}_{ij}}\circ e^{\tilde{X}_i}=e^{\rho(\rho(\rho(-\tilde{X}_j,T_{ij}),\tilde{X}_{ij}),\tilde{X}_i)}
\end{equation*}
The relation (\ref{ro2}) provides the following expansion:
\begin{multline}\label{rel3}
\rho(\rho(\rho(-\tilde{X}_j,T_{ij}),\tilde{X}_{ij}),\tilde{X}_i)=\\ T_{ij}+\tilde{T}_{ij}+\underbrace{\fraction{1}{2}\left[T_{ij},\tilde{X}_j\right]+\fraction{1}{2}\left[T_{ij},\tilde{X}_{ij}\right]+\fraction{1}{2}\left[T_{ij},\tilde{X}_i\right]+\cdots}_{Y_{ij}}.
\end{multline}
Now, for any integers $n$ and $m$,
$ \left[\mathfrak{M}^{m}_Z\mathfrak{X}_{S,Z},\mathfrak{M}^{n}_Z\mathfrak{X}_{S,Z}\right]\subset\mathfrak{M}^{m+n}_Z\mathfrak{X}_{S,Z}
$. Hence, $\coup{Y_{ij}}$ is a $1$-cocycle of
$\ZZ^1\left(\UU,\mathfrak{M}^{n+1}_Z\mathfrak{X}_{S,Z}\right)$. So, we get the following expression
\begin{eqnarray}
e^{-\tilde{X}_j}\circ\phi_j^{-1}\circ\phi_{{ij}}\circ \phi_i\circ e^{\tilde{X}_i}&=&e^{T_{ij}+\tilde{T}_{ij}+Y_{ij}}\circ \phi_{ij}^1\circ\left[\phi_{ij}^1,e^{-\tilde{X}_i}\right]\\
&=&e^{T_{ij}+\tilde{T}_{ij}}\circ\phi_{ij}^2\circ \phi_{ij}^1\circ \left[\phi_{ij}^1,e^{-\tilde{X}_i}\right].\label{rel4}
\end{eqnarray}
where $\phi_{ij}^2$ is equal to $e^{-T_{ij}-\tilde{T}_{ij}}\circ e^{T_{ij}+\tilde{T}_{ij}+Y_{ij}}$. Again, the Campbell-Hausdorff formula shows that $\phi^2_{ij}$ is the flow of the vector field admitting the following expansion
$Y_{ij}-\frac{1}{2}\left[Y_{ij},T_{ij}+\tilde{T}_{ij}\right]+\cdots$. As a consequence, in view of the property (\ref{cle}), $\phi^2_{ij}$ takes its values in $\GG^{n+1}_Z$. Finally, similary arguments ensure that the commutator in (\ref{rel4}) is in $\GG^{n+1}_Z$ class.
Hence, the equation (\ref{rel4}) is the induction hypothesis at rank $n+1$. 

\end{demo}
Let us now prove the proposition (\ref{Glulg1}). One can write
$ \NN\simeq {\MM}[\phi_{ij}]$ with $\coup{\phi_{ij}}$ in
$\ZZ^1(\UU,\GG_Z^1)$.  The infinitesimal cobordism lemma gives us three cocycles $\coup{T_{ij}}$ in $\ZZ^1\left(\UU,\mathfrak{M}_Z\mathfrak{X}_{\FF,Z}\right)$,
$\coup{\phi_i}$ in $\ZZ^0\left(\UU,\GG_Z^1\right)$ and $\coup{\tilde{\phi}_{ij}}$ in $\ZZ^1(\UU,\GG_Z^{n})$ such that
$
\phi_j^{-1}\phi_{ij}\phi_i=e^{T_{ij}}\tilde{\phi}_{ij}
$.
Now, the (\ref{lell}) and (\ref{lell1}) properties imply that
$$ \MM[\phi_{ij}]\stackrel{\mathtt{G}_1}{\simeq} {\MM}[\phi_j^{-1}\phi_{ij}\phi_i]= {\MM}[e^{T_{ij}}\tilde{\phi}_{ij}]\stackrel{\mathtt{G}_1}{\simeq}{\MM}[e^{T_{ij}}][\zeta^1\tilde{\phi}_{ij}].$$
Moreover, the stability property applied in $\text{Glu}_n(\MM[e^{T_{ij}}],Z[e^{T_{ij}}],\UU[e^{T_{ij}}])$ to the cocycle $\coup{\zeta^1\tilde{\phi}_{ij}}$ shows that
$ {\MM}[e^{T_{ij}}][\zeta^1\tilde{\phi}_{ij}]\stackrel{\mathtt{G}_1}{\simeq} {\MM}[e^{T_{ij}}]$.
Hence,
$ {\MM}[\phi_{ij}]$ is isomorphic to ${\MM}[e^{T_{ij}}]$,
which is the checked property.

\subsubsection{Step 3:  cobordism in $\textup{Glu}_0(\MM,Z,\UU)$.}\label{induction}

The proof of the (\ref{cobordisme0}) proposition is an induction on the height of the tree, which allows us to use the previous result established for the $\textup{Glu}_1(\MM,Z,\UU)$ category.

\begin{figure}[ht]\label{ind}

\begin{center}
\includegraphics{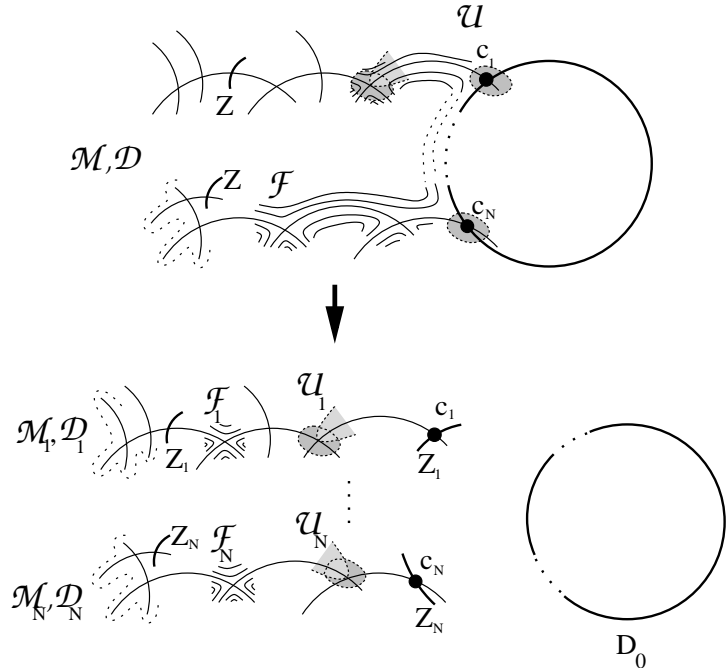}
\end{center}{\bf \caption{Induction construction.}}
\end{figure}

\medskip
\noindent Let $\coup{\phi_{ij}}$ be such that $\NN=\MM[\phi_{ij}]$.
Let us denote by $D_0$ the irreducible component of $\DD$ appearing after the first blowing-up. Let
$\fami{c_1,\ldots,c_N}$ be the set of singular points of $\DD$ on
$D_0$. We denote by $\DD_l$ the branch of $\DD$ attached to $c_l$ and $\MM_l$ the neighborhood of $\DD_l$ in $\MM$. $\UU_l$ refers to the distinguished covering of $\DD_l$ induced by restriction of $\UU$. Moreover, we denote by $U_0\in \UU$ the open set $D_0\backslash\fami{c_1,\ldots,c_N}$ and $U_l\in \UU$ the neighborhood of $c_l$.
The foliation $\FF_l$ refers to the restriction of the foliation
$\FF$ to $\MM_l$. We consider a cross $Z_l$ 
defined by the trace in $\MM_l$ of $D_0$ and of the strict transform of $Z$ (see Fig. \ref{ind}). 
One can see that $\FF_l$ is a non-dicritical generalized curve on a marked crossed tree and that $Z_l$ is adapted to $\FF_l$.

\noindent Let us consider $\coup{\phi_{ij}^l}\in \ZZ^1(\UU_l,\GG^0_{Z_l})$ the $1$-cocycle restriction of $\coup{\phi_{ij}}\in \ZZ^1(\UU,\GG^0_{Z})$ to $\UU_l$,
$$ \phi_{ij}^l=\phi_{ij},\quad U_{i}, U_j\in \UU_l .$$
Since each tree $\DD_l$ has a height smaller than $h-1$, the induction hypothesis ensures that each tree
$\MM_l[\phi_{ij}^l]$ is $\FF_l$-cobordant to $\MM_l$ in the category $\textup{Glu}_0(\MM_l,Z_l,\UU_l)$.

\noindent Let us assume first that the above $\FF_l$-cobordisms are elementary. By definition, there exists a family of $1$-cocyles $\coup{T^l_{{ij}}}$, $l=1,\ldots,N$ of $\ZZ^1(\UU_l,\mathfrak{I}_{Z_l}\mathfrak{X}_{\FF_l,Z_l})$  such that,
\begin{equation*}
 \MM_l[\phi^l_{ij}]\stackrel{\mathtt{G}_0}\simeq  \MM_l[e^{T^l_{ij}}].
\end{equation*}
In view of the lemma (\ref{lell}), there exists a family of $0$-cocycles in $\GG^0_{Z_l}$ verifying
\begin{equation}\label{eq}
\phi^l_{ij}=\psi^l_je^{T^l_{ij}}{\psi^l_i}^{-1}.
\end{equation}
Let us suppose that $\psi^l_1$ refers to the component of the $0$-cocycle $\coup{\psi^l_i}$ defined on the open set $U_l$ containing the singularity $c_l$. Let $\coup{e^{T_{ij}}}_{ext}$ be the
$1$-cocyle 
$$
\coup{e^{T_{ij}}}_{ext}= \left\{\begin{array}{lr}
    \text{Id} & \text{on } U_l\cap U_0 \\
    e^{T^l_{{ij}}} & \text{else}
\end{array}\right..$$ We denote by $\phi_{l0}$ the component of the $1$-cocycle $\coup{\phi_{{ij}}}$ defined on $U_l\cap U_0$. Finally, we define a $1$-cocycle by
$$
\kappa_{{ij}} = \left\{\begin{array}{lr}
    \psi^l_1\phi_{l0} & \text{on } U_l\cap U_0 \\
    e^{T^l_{{ij}}} & \text{else}
\end{array}\right. . $$
The relation (\ref{eq}) induces an isomorphism in $\text{Glu}_0(\MM,Z,\UU)$
\begin{equation}\label{Gluiso}
\MM[\kappa_{{ij}}]\stackrel{\mathtt{G}_0}\simeq \MM[\phi_{{ij}}].
\end{equation}
We define the following $1$-cocycle 
$$
\tilde{\kappa}_{{ij}} = \left\{\begin{array}{lr}
\psi^l_1\phi_{l0} & \text{on } U_l\cap U_0 \\
\text{Id} & \text{else}
\end{array}\right. .$$
Let $\tilde{\MM}\in\text{Glu}_0(\MM,Z,\UU)$ be defined by $\tilde{\MM}=\MM[e^{T_{ij}}_{ext}]$. We have

\begin{eqnarray}
\tilde{\MM}[\zeta^1\tilde\kappa_{{ij}}] & = &  \MM[e^{T_{{ij}}}_{\text{ext}}][\zeta^1\tilde\kappa_{{ij}}]\\
& \stackrel{\mathtt{G}_0}\simeq &  \MM[e^{T_{ij}}_\text{ext}\tilde\kappa_{ij}]\\
& = & \MM[\kappa_{{ij}}]\label{egaliteGlu}.
\end{eqnarray}
Outside $U_l\cap U_0$, the components of the cocycle $\coup{\zeta^1\tilde\kappa_{{ij}}}$ are equal to $\text{Id}$. Since $\psi_l$ is in the 
$\GG^0_{Z_l}$ class, $\zeta^1\tilde\kappa_{{l0}}$ is a germ of automorphism in $\GG^0_{\tilde{Z}}$ where $\tilde{Z}=Z[e^{T_{ij}}_\text{ext}]$. Now, the sheaves
$\mathfrak{I}_{\tilde{Z}}$ and $\mathfrak{M}^1_{\tilde{Z}}$ coincide along the regular part of $\tilde{D}_0$. Hence, the sheaves 
$\GG^1_{\tilde{Z}}$ and $\GG^0_{\tilde{Z}}$ are equal too. As a consequence, the $1$-cocycle
$\coup{\zeta^1\tilde\kappa_{ij}}$ is in $\GG_{\tilde{Z}}^1$ and $\tilde{\MM}[\zeta^1\tilde\kappa_{ij}]$ belongs to $\text{Glu}_1(\tilde{\MM},\tilde{Z},\tilde{\UU})$. The proposition
(\ref{Glulg1}) ensures that there exists an elementary $\FF[e^{T_{ij}}_\text{ext}]$-cobordism between $\tilde{\MM}$ and $\tilde{\MM}[\zeta^1\tilde\kappa_{ij}]$ defined by a $1$-cocycle $\coup{e^{\tilde{T}_{ij}}}\in Z^1(\tilde{\UU},\mathfrak{I}_{\tilde{Z}}\XX_{\FF[e^{T_{ij}}_\text{ext}],\tilde{Z}})$. Thanks to (\ref{Gluiso}) and (\ref{egaliteGlu}), we have $\NN\stackrel{\mathtt{G}_0}{\simeq}\tilde{\MM}[\zeta^1\tilde\kappa_{ij}]$. Hence, we get a final cobordism between $\NN$ and $\MM$:
\begin{eqnarray*}
\NN\stackrel{\mathtt{G}_0}{\simeq}\tilde{\MM}[e^{\tilde{T}_{ij}}], \\
\tilde{\MM}\stackrel{\mathtt{G}_0}{\simeq} \MM[e^{T_{ij}}_\text{ext}].
\end{eqnarray*}
\noindent Now, if the cobordisms are not elementary, one can suppose them to be decomposed in sequences of elementary cobordisms of same length. Indeed, we have
\begin{eqnarray*}\label{long}
 \MM_1&=&\MM^1_1\fleche{\FF_1^1,Z_1^1}\MM_1^2\cdots \MM_1^{p-1} \fleche{ \FF^{p-1}_1,Z^{p-1}_1}\MM^{p}_1\stackrel{\mathtt{G}_0}{\simeq} \NN_1,\\
  \MM_2&=&\MM^1_2\fleche{\FF_2^1,Z_2^1}\MM_1^2\cdots \MM_2^{p-1} \fleche{ \FF^{p-1}_2,Z^{p-1}_2}\MM^{p}_2\stackrel{\mathtt{G}_0}{\simeq} \NN_2,\\
&\vdots & \\
 \MM_N&=&\MM^1_N\fleche{\FF_N1,Z_N^1}\MM_N^2\cdots \MM_N^{p-1} \fleche{ \FF^{p-1}_N,Z^{p-1}_N}\MM^{p}_N\stackrel{\mathtt{G}_0}{\simeq} \NN_N,
\end{eqnarray*}
where $\FF^k_l$ and $Z^k_l$, $l\geq 1$ naturally refer to the foliations and the crosses induced by the successive cobordisms. Let $\coup{e^{T^{k,l}_{ij}}}$ be the cocycle defining the elementary cobordism $\MM^k_l\fleche{\FF^k_l,Z^k_l} \MM_{l}^{k+1}$, where  $\coup{T^{k,l}_{ij}}$ belongs to $\ZZ^1(\UU^k_l,\mathfrak{I}_{Z^k_l}\XX_{\FF^k_l,Z^k_l})$.
For $1\leq k\leq p-2$, we consider the trees $\dot{\MM}^k$ defined by $\dot{\MM}^1=\MM$ and
 $$\dot{\MM}^{k+1}=\dot{\MM}^k[{e^{T^{k,l}_{ij}}}_{ext}]$$ with $$
\coup{e^{T^{k,l}_{ij}}}_{ext}= \left\{\begin{array}{lr}
    \text{Id} & \text{on } U^{k,l}_l\cap U^{k,l}_0 \\
    e^{T^{k,l}_{{ij}}} & \text{else}
\end{array}\right. .$$
Here $\UU^{k,l}$ refers to the distinguished covering induced by the successive gluings, $\UU^{k+1,l}=\UU^{k,l}[{e^{T^{k,l}_{ij}}}_{ext}]$. In view of the construction, the previous relations give us a general 
$\FF$-cobordism between $\MM$ and $\dot{\MM}^{p-1}$. Moreover, for all $l$ the tree $\dot{\MM}^{p-1}_l$, restriction of the tree $\dot{\MM}^{p-1}$ over
 the singularity $l$ is isomorphic to $\MM^{p-1}_l$ and the foliation $\dot{\FF}^{p-1}_l$ got by successive cobordisms is isomorphic to $\FF^{p-1}_l$. Hence, to finish the proof of (\ref{cobordisme0}), one must solve the elementary case from $\dot{\MM}^{p-1}$ to $\NN$, what has already been done.

\subsection{The second kind case.}\label{deuxspe}
In this subsection, we want to extend the proposition (\ref{cobordisme0}) to a bigger and more natural class of foliations than the generalized curves. 

\medskip
\noindent A reduced singularity of foliation with a linear part of the form $x\dd y$ admits formal normal forms of following types \cite{diff}:
$$x^{p+1}\dd y-y(1+\lambda x^p)\dd x.$$
The separatrix ${x=0}$ is called the strong invariant curve. If some singularities of such kind appear in the reduction process, most of previous results are false, 
except when this kind of singularities is in a specific position.
The notion of foliation of second kind was introduced in \cite{MatSal} to take care of that possibility. This class of foliations admits the same properties as generalized curves and provides a nice description of local formal moduli spaces. 
\begin{definition}
$\FF$ is of second kind if
\begin{enumerate}
\item $\FF$ is non-dicritical,
\item each singularity of the divisor is a reduced singularity of the foliation with two non-vanishing eigenvalues,
\item each singularity of the foliation on the regular part of the divisor is either of above kind or has $x\dd y$  for linear part; in the latter case, one asks the strong invariant curve to be the germ defined by the divisor. 
\end{enumerate}
\end{definition}
\noindent With obvious notations, we consider the sheaf $\widehat{\mathfrak{X}}_{S,Z}$ over $\DD$ of germs of formal vector field that are tangent to the divisor, the separatrix and the cross. We denote 
$\widehat{\mathfrak{X}}_{\FF,Z}$ the subsheaf of vector fields tangent to the foliation and the cross. In \cite{MatQuasi}, one can find the following criterion for a foliation to be of second kind:

\begin{prop}
The foliation $\FF$ is of second kind if and only if the sequence of sheaves
$$0 \longrightarrow \widehat{\mathfrak{X}}_{\FF,Z} \longrightarrow \widehat{\mathfrak{X}}_{S,Z}\xrightarrow{E^*\omega(.)} \widehat{\mathcal{O}} \left(f\circ E\right)\longrightarrow 0 $$
is exact.
\end{prop}
\noindent Now, it's not hard to see that from the above exact sequence, every previous constructions and arguments can be repeated in the formal context under the second kind hypothesis. More precisely, one can get the following result:
\begin{prop}\label{cobordismeformel}
Let $\FF$ be a foliation of second kind on a marked tree $\MM$ and $Z$ an adapted cross. Let $\NN$ be in $\textup{Glu}_0(\MM,Z,\UU)$. There exists a finite sequence of $1$-cocyles
$\coup{\widehat{T}^k_{ij}}_{k=1..N}$ such that
$$\NN\stackrel{\widehat{\mathtt{G}_0}}{\simeq}\MM[e^{\widehat{T}^1_{ij}}][e^{\widehat{T}^2_{ij}}][\ldots][e^{\widehat{T}^N_{ij}}].\footnote{Here, the isomorphism is related to the category $\widehat{\textup{Glu}}_0(\MM,Z,\UU)$. This notation refers to the same one but in convergent context. The transposition to formal context is easy.}$$
Here, $\coup{\widehat{T}^p_{ij}}$ is a $1$-cocycle with values in
$\widehat{\mathfrak{X}}_{\widehat{\FF}_p,\widehat{Z}_p}$: this is the sheaf over the divisor of the tree
$\MM[e^{\widehat{T}^1_{ij}}][e^{\widehat{T}^2_{ij}}][\ldots][e^{\widehat{T}^{p-1}_{ij}}]$ whose fiber is the space of germ of formal vector field tangent to the following formal foliation and formal cross
$$ \widehat{\FF}_p=\widehat{\FF}[e^{\widehat{T}^1_{ij}}][e^{\widehat{T}^2_{ij}}][\ldots][e^{\widehat{T}^{p-1}_{ij}}],
 \quad \widehat{Z}_p=Z[e^{\widehat{T}^1_{ij}}][e^{\widehat{T}^2_{ij}}][\ldots][e^{\widehat{T}^{p-1}_{ij}}].$$
\end{prop}
\noindent From this formal construction, one can go back to the convergent context using the following lemma: let $\MM$ be foliated by a convergent foliation $\FF$ of second kind and let $\coup{\widehat{T}_{ij}}$ be in $Z^1(\UU,\widehat{\mathfrak{I}}_Z\widehat{\mathfrak{X}}_{\FF,Z})$.
\begin{lemme}\label{depfor}
There exists  $\coup{T^c_{ij}}\in Z^1(\UU,\mathfrak{I}_Z\mathfrak{X}_{\FF,Z})$ such that
$\MM[\widehat{T}_{ij}]\stackrel{\widehat{\mathtt{G}_0}}{\simeq}\MM[T^c_{ij}]$.
\end{lemme}
\begin{demo} Since $\FF$ is convergent, there exists a convergent vector field $T_{ij}$ tangent to the foliation and a formal series $\widehat{\phi}_{ij}$
such that
$\widehat{T}_{ij}=\widehat{\phi}_{ij}T_{ij}$.
For any integer $n$, we consider the cocyle $\coup{T^c_{ij}}$ in $\mathfrak{X}_{\FF,Z}$ defined by
$T^c_{ij}=\phi^n_{ij}T_{ij}$,
where $\phi^n_{ij}$ refers to a representative function of the $n$-jet of $\widehat{\phi}_{ij}$. We are going to show that, for $n$ big enough, the gluings associated to $\widehat{T}_{ij}$ and $T^c_{ij}$ are formally equivalent in $\widehat{\textup{Glu}}_0(\MM,Z,\UU)$. Let us consider the following $1$-cocycle
$$\delta_{ij}(x,t,s)=\big(e^{(t-s)\widehat{T}_{ij}}\circ
e^{(s)T^c_{ij}},t,s\big).$$  We construct an equisingular unfolding in the following way:

\begin{equation*}
\MM_{\overline{\mathbb{D}}^2}=\bigcup_i
U_i\times\overline{\mathbb{D}}^2\left/_{(x,t,s)\sim \delta_{ij}(x,t,s)}\right.,
\end{equation*}
which admits a projection 
$\MM_{\overline{\mathbb{D}}^2}\xrightarrow{\Pi}\overline{\mathbb{D}}^2$.
This manifold is the neighborhood of a divisor $\DD_{\overline{\mathbb{D}}^2}$. Thanks to local triviality of the unfolding along $U_i\times\overline{\mathbb{D}}^2$, one can find a family $\fami{X_i}$ of vector fields tangent to the foliation such that in
$U_i\times\overline{\mathbb{D}}^2$ we have
$T\Pi(X_i)=\fraction{\partial}{\partial s}$.
The $1$-cocycle 
$\coup{X_i-X_j}$
takes its values in the sheaf over $\DD_{\dbar}$ of
germs of vector fields tangent to the foliation, vertical - $T\Pi(X)=0$ - and vanishing at order $n$ along $\DD_{\dbar}$. We denote by $\widehat{\mathfrak{X}}_{\FF,Z,n}^{\overline{\mathbb{D}}^2}$  the latter sheaf. Since $\widehat{\mathfrak{X}}^{\dbar^2}_{\FF,Z,1}$ is coherent, one can find the following property in \cite{camacho2} :

\begin{souslemme}
For $n$ big enough, the map
$$ H^1\left(\widehat{\mathfrak{X}}_{\FF,Z,n}^{\overline{\mathbb{D}}^2}\right)\longrightarrow
  H^1\left(\widehat{\mathfrak{X}}_{\FF,Z,1}^{\overline{\mathbb{D}}^2}\right)$$
is trivial.
\end{souslemme}
\noindent Hence, there exists $0$-cocycle $\coup{Y_i}$ in
$\widehat{\mathfrak{X}}_{\FF,Z,1}^{\overline{\mathbb{D}}^2}$ such that
$$X_i-X_j=Y_i-Y_j$$
Therefore, $X=X_i-Y_i$ is a global formal vector field such that
$T\Pi(X)=\fraction{\partial}{\partial s}$.
The biholomorphim
$\phi(x,t,s)=(e^{(s)X}(x,t,0),t,s)$ formally conjuguates the equisingular unfoldings
$\MM_{\overline{\mathbb{D}}^2}$ and
$\MM_{\overline{\mathbb{D}}\times\{0\}}\times\overline{\mathbb{D}}$. By restricting them along the diagonal, one can see that 
$\MM[\widehat{T}_{ij}]$ and $\MM[T^c_{ij}]$ are conjugated.
\end{demo}
By applying the previous lemma to each elementary cobordism in (\ref{cobordismeformel}), one can proof the cobordism result for foliation of second kind: 
\begin{prop}\label{cobordisme0formel}
Let $\FF$ be a foliation of second kind on a marked tree $\MM$ and $Z$ a cross adapted to $\FF$. Any
element in $\textup{Glu}_0(\MM,Z,\UU)$ is $\FF$-cobordant to $\MM$.
\end{prop}

\section{Existence theorem.}\label{real}
In this section, we use the proposition (\ref{cobordisme0formel}) to establish the existence theorem. Let $\FF$ be a foliation of second kind on a marked tree $(\MM,\DD,\Sigma_\FF)$ and  $(\MM',\DD',\Sigma')$ be any marked tree.
\begin{thm}[Existence theorem]\label{thmA}
\

\noindent If the marked weighted dual trees of $\MM$ and $\MM'$ are conjugated then there exists a foliation $\FF'$ on  $\MM'$ 
such that $\FF$ and
$\FF'$ are linked by an equisingular unfolding, which respects the marking.
\end{thm}
\noindent We focus on the proof of the above statement.

\subsection{First step: cocyle transformation.}
Since the marked dual trees of $\MM$ and $\MM'$ are conjugated, one can find in \cite{Seguy} the following result: there exists an equisingular unfolding of $\FF$, which leads to a foliation defined on a tree with a divisor biholomorphic to $\MM'$'s one. Hence, one can suppose the divisors $\DD$ and $\DD'$ isomorphic. Let us denote by $h$ a biholomorphism between the divisors. For each component $D$, we consider a fibration $\pi_{D}$ transversal to $D$
$$\pi_{D}:T(D)\longrightarrow D$$ where $T(D)$ is a fixed tubular neighborhood of $D$. Futhermore, we assume the trace of transversal components of $\DD$ to be some fibers of
$\pi_{D}$. We make the same construction over
$\DD'$'s components. An easy computation in coordinates allows us to show the following lemma:
\begin{lemme}
There  exists a family of maps $\coup{H_{D}}_{D\in Comp(\DD)}$ such that the following diagram is commutative
\begin{center}
$\begin{CD}
  T(D) @>H_{D}>> h \left ( T(D) \right) \\
  @V\pi_{D}VV @V\pi_{h(D)}VV\\
  D @>h_{|D}>> h(D)
\end{CD}
$
\end{center}
\end{lemme}

\noindent We denote by $\text{Comp}(\DD)^{\check{2}}$ the set $\fami{(D,D')\in \text{Comp}(D)^2 | D\cap D'\neq \emptyset}$. Let us consider the germ of manifold defined by 
$$\tilde{\MM}=\left. \coprod_{D\in Comp(\DD)} T(D)\right/ \coup{x\sim {H_{D}}^{-1}\circ H_{D'}x}_{(D,D')\in\text{Comp}(\DD)^{\check{2}} }.$$
One has to notice that the above gluing is different from the gluing introduced in previous sections: here, coverings are made of tubular neighborhoods, which intersect each other along polydisks; whereas open sets of distinguished covering are finer and intersect each other along torus. Nevertheless, the tree $\tilde{\MM}$ is a neighborhood of some divisor $\tilde{\DD}$. The tree $(\tilde{\MM},\tilde{\DD},\tilde{\Sigma})$ is a marked tree conjugated to $(\MM',\DD',\Sigma')$: indeed the family $\coup{H_{D}}$ induces a biholomorphism adapted to markings.

\noindent In order to apply the proposition (\ref{cobordisme0}), we are going to build a tree $\dot{\MM}$ such that $\dot{\MM}$ verifies the existence theorem and $\tilde{\MM}$ belongs to $\text{Glu}_0(\dot{\MM},Z,\UU)$ for suitable cross and covering.

\noindent Let us denote by $s_{DD'}$ the intersection of $D$ and $D'$ and $\phi_{DD'}:={H_{D}}^{-1}\circ H_{D'}$.
\begin{lemme}\label{lelele}
There exists families $\coup{\Delta_{DD'}}_{(D,D')\in \textup{Comp}(\DD)^{\check{2}}}$ and 
$\coup{\phi_D}_{D\in \textup{Comp}(\DD)}$ of automorphism germs such that
\begin{itemize}
\item $\Delta_{DD'}$ is defined near $s_{DD'}$, fixes this point and lets invariant each local leaf; 
\item $\phi_D$ is defined along $D$ and fixes each point of $D$;
\item the germ of automorphism $\phi_{D'}^{-1}\circ\Delta_{DD'}\circ\phi_{DD'}\circ \phi_D$
is tangent to the identity at $s_{DD'}$.
\end{itemize}
\end{lemme}
\begin{demo}  Let us consider the standard metric $d$ on the dual weighted graph of $\MM$ and fix one vertex $D_0\in\text{Comp}(\DD)$. We define the subgraph $\mathbb{A}^*_n, n\geq 0$ whose vertex are at distance smaller than $n$ from $D_0$. The graphs $\mathbb{A}^*_n$ are connected and cover the whole $\MM$'s graph. One can now consider the restricted family $\coup{\phi^{n}_{DD'}}$ defined by 
$$\phi^{n}_{DD'}=\phi_{DD'},\quad DD' \text{ is an edge of } \mathbb{A}^*_n$$
By induction, we establish the result on the subfamily $\coup{\phi^n_{DD'}}$. For $n=0$, the result is obvious and one can choose for $\phi_{D_0}$ the identity automorphism. Let us suppose the result true for $n$: precisely, we have two families $\Delta_{DD'}$ and $\phi_{D}$ as in the lemma such that for any edge $DD'$ of $\mathbb{A}^*_n$ 
$$\phi_{D'}^{-1}\circ\Delta_{DD'}\circ\phi_{DD'}\circ \phi_D$$ is tangent to the indentity at $s_{DD'}$. Let $D_iD_j$ be an edge of $\mathbb{A}^*_{n+1}\backslash\mathbb{A}^*_{n}$ such that $D_i$ is a vertex of $\mathbb{A}^*_n$ and $D_j$ a vertex of $\mathbb{A}^*_{n+1}$. One has to define in a good way $\phi_{D_j}$ and $\Delta_{D_iD_j}$ in order to conclude. Now, one can easily prove the general following results:
\begin{souslemme}\label{transfo1}\begin{enumerate}
\item
Let $s$ be a singular point of  $\DD$ and $\alpha \in \mathbb{C}^*$. There exists
a germ of automorphism $\Delta$ defined near $s$ and letting invariant each local leaf such that, in adapted coordinates, the tangent map at $s$ is 
$$T_s\Delta =\left(\begin{array}{cc} \alpha & 0\\ 0 & * \end{array}\right).$$
\item 
For any $\beta$ in $\mathbb{C}^*$, any component of $D$ and any $c\in D$, there exists a germ of automorphism $\phi$ defined in the neighborhood of $D$ fixing each point of $D$ such that, in adapted coordinates\footnote{Here, we assume the second coordinates to be transversal to the divisor}, the tangent map at $c$ is
$$T_c\phi =\left(\begin{array}{cc} 1 & 0\\ 0 & \beta \end{array}\right).$$
\end{enumerate}
\end{souslemme}
\noindent In adapted coordinates, the tangent map of the automorphism $\phi_{D_iD_ j}\circ \phi_{D_i}$ is
$\left(\begin{array}{cc} \alpha & 0\\ 0 & \beta \end{array}\right)$. In view of the previous lemma,
there exists a germ of automorphism $\Delta_{D_iD_j}$ defined near $s_{D_iD_j}$ letting invariant each local leaf and a germ of automorphism $\phi_{D_j}$ defined in the neighborhood of $D_{j}$ fixing each point of $D_{j}$ such that
$$T_{s_{D_iD_j}}(\Delta_{D_iD_j}^{-1}\circ \phi_{D_j})=\left(\begin{array}{cc} \alpha & 0\\ 0 & \beta \end{array}\right).$$
Then, 
$T_{s_{D_iD_j}}( \phi_{D_j}^{-1}\circ\Delta_{D_iD_j}\circ\phi_{D_iD_j}\circ \phi_{D_i})=\text{Id}$. Since the dual tree doesn't have any cycle, one can repeat the same construction for all edges of $\mathbb{A}^*_{n+1}\backslash\mathbb{A}^*_{n}$. Hence, the lemma is proved for $\mathbb{A}^*_{n+1}$.
\end{demo}

\subsection{Second step: a fine distinguished covering.}
With the family $\fami{\Delta_{DD'}}$ of the previous lemma, one can define
$$\dot{\MM}=\left. \coprod_{D\in Comp(\DD)} T(D)\right/  \coup{x\sim\Delta^{-1}_{DD'}x}.$$ The tree $\dot{\MM}$ is a neighborhood of some divisor $\dot\DD$. One can find in \cite{iso} a precise description of germs letting invariant each local leaf of a reduced singularity with two non-vanishing eigenvalues. These can be written $(x,y)\mapsto \phi^{t(x,y)}_X$ where $t(x,y)$ is a function, $X$  a germ of vector field tangent to the foliation and $\phi^{t}_X$ the flow of $X$ at time $t$. Particulary, this description ensures that $\dot{\MM}$ admits a foliation $\dot{\FF}$ linked to $\FF$ by an equisingular unfolding. Basically, the unfolding exists because there is a suitable isotopy from $(x,y)\mapsto \phi^{t(x,y)}_X$ to the identity defined by $(s,x,y)\mapsto \phi^{st(x,y)}_X$. 

\noindent One notices that we have the following isomorphism
$$ \tilde{\MM}\simeq\left. \coprod_{\dot{D}\in Comp(\dot\DD)} T(\dot{D})\right/\coup{x\sim\Delta_{DD'}\circ\phi_{DD'}x}.$$
Let us consider $\theta_{DD'}:={\phi_{D'}}^{-1}\circ\Delta_{DD'}\circ\phi_{DD'}\circ\phi_D$. One can see that, 
$$\left. \coprod_{\dot{D}\in Comp(\dot\DD)} T(\dot{D})\right/ \coup{x\sim\theta_{DD'}x}\simeq \tilde{\MM}.$$

\noindent We are going to make a last transformation on $\theta_{DD'}$ in order to obtain a cocycle with values in $\GG^0_Z$; it's not hard to prove the following result:
\begin{souslemme}
Let $\theta$ be a germ of automorphism tangent to identity at
$0\in\mathbb{C}^2$ letting invariant the axes $\{x=0\}$ and $\{y=0\}$. $\theta$ admits a decomposition of the form
$$\theta= \theta^0\circ \theta^1, \quad  \theta^0=\left(\begin{array}{c} x+x^2(cldots) \\ y+xy(\cdots) \end{array}\right),\quad  \theta^1= \left(\begin{array}{c} x+xy(\cdots) \\ y+y^2(cldots) \end{array}\right).$$
\end{souslemme}
\noindent For each $D,D'\in \text{Comp}(\DD)^{\check{2}}$,
we decompose $\theta_{DD'}$ in $\theta_{DD'}=\theta_{DD'}^0\circ\theta_{DD'}^1$.
Taking a distinguished covering $\UU=\coup{U_i}$ finer than the tubular covering $\coup{T(\dot{D})}_{D\in \text{Comp}(\DD)}$, we have the following isomorphisms 
$$\tilde{\MM}\simeq\left.\coprod T(\dot{D})\right/\theta_{DD'}\simeq
\left.\coprod U_i\right/\theta_{DD'}^\epsilon=\dot{\MM}[\theta^\epsilon_{DD'}].$$
Let $Z$ be any suitable cross on $\dot{\MM}$ adapted to $\dot{\FF}$. Clearly, the tree $\tilde{\MM}$ belongs to $\text{Glu}_0(\dot{\MM},Z,\UU)$. In view of the previous results, $\tilde{\MM}$ is $\dot{\FF}$-cobordant to $\dot{\MM}$. Hence, $\MM'\simeq \tilde{\MM}$ admits a foliation $\FF'$ linked to $\FF$ by an equisingular unfolding. 
\carre

\section{Cobordism class.}\label{princp}
Let $\FF$ be a non-dicritical generalized curve and $S=\text{Sep}(\FF)$ its separatrix. A foliation $\FF'$ is said to be cobordant to $\FF$ when there is an equisingular unfolding in the sense of \cite{univ} that links $\FF$ and $\FF'$. We denote by $\textup{Cob}(\FF)$ the set of cobordant foliation. The set $\text{E}(S)$ refers to the equisingularity class of $S$ or, in a equivalent way, the set of curves topologically conjugated to $S$.
\begin{thm}\label{thmB}
 The map defined by
$$\FF'\in \textup{Cob}(\FF)\longmapsto \textup{Sep}(\FF')\in\textup{E}(S)$$
is onto.
\end{thm}
\noindent In order to prove the above result, one notices that the existence theorem holds for trees that may not be the minimal reduction tree of a foliation. This remark leads us to define the following sequence of blowing-up process: let  $\coup{E^j}_{j=1\ldots h}$ be the blowing-up process of the reduction of the foliation $\FF$ (\ref{ArMr}). We define a sequence of blowing-up process in this way:
\begin{itemize}
\item $\mathfrak{E}_h=\coup{E^j}_{j=1\ldots h}$.
\item $\mathfrak{E}_{n+1}$ is the blowing-up process build over $\mathfrak{E}_n$ where $S_n$ and $\Sigma_n$ are both the set of singularities of $E_n^*\FF$, where $E_n$ refers to the total morphism of $\mathfrak{E}_n$ (\ref{ArMr}).
\end{itemize}

\begin{demo}{\it\  theorem (\ref{thmB})} Let $S'$ be a curve in the equinsigularity class of $S$.
Let $E$ and $E'$ refer respectively to the reduction process of $S$ and $S'$. In view of a result in \cite{camacho}, since $\FF$ is a generalized curve, $E$ is also the reduction of $\FF$. Moreover, according to a result of O. Zariski \cite{zariski}, the weighted dual graphs of $S$ and $S'$ are conjugated. Particulary, the attached points of $S$ and $S'$ strict transforms are on conjugated components.

\noindent 
Let $\mathcal{H}$ be the foliation $\dd h=0$ where $h$ is any reduced equation of $S'$. For any integer $n$ bigger than the height of $E$, we consider the sequence $E_n$ and $E'_n $ built as above with $\FF$ and $\mathcal{H}$ as respective initial data. For every $n$, the weighted dual trees associated to $E_n$ and $E'_n$ are conjugated. The theorem (\ref{thmA}) ensures the existence of a foliation $\FF_n$ on the tree of $E'_n$ linked to $E_n^*\FF$ by an equisingular unfolding. Let $S_n$ be the germ of separatrix of $\FF_n$ at $0\in \mathbb{C}^2$. One can see that the attached points on the divisor of the strict transforms by $E'_n$ of the curves $S_n$ and $S'$ are on same components. Since $\FF_n$ is cobordant to $\FF$, $\FF_n$ is topologically equivalent to $\FF$. In particular, $S_n$ is topologically equivalent to $S'$. For $n$ great enough, the study in \cite{zariski} ensures that $S_n$ and $S'$ are analytically conjugated. The image of $\FF_n$ by this conjugacy is a holomorphic foliation topologically equivalent to $\FF$ admitting $S'$ for separatrix.
\end{demo}

\bibliography{biblio}

\def\cprime{$'$} \def\cprime{$'$}
\begin{thebibliography}{10}

\bibitem{iso}
M.~Berthier, D.~Cerveau, and R.~Meziani.
\newblock Transformations isotropes des germes de feuilletages holomorphes.
\newblock {\em J. Math. Pures Appl. (9)}, 78(7):701--722, 1999.

\bibitem{camacho}
C.~Camacho, A.~Lins~Neto, and P.~Sad.
\newblock Topological invariants and equidesingularization for holomorphic
  vector fields.
\newblock {\em J. Differential Geom.}, 20(1):143--174, 1984.

\bibitem{camacho2}
C.~Camacho and H.~Movasati.
\newblock {\em Neighborhoods of analytic varieties}, volume~35 of {\em
  Monograf\'ias del Instituto de Matem\'atica y Ciencias Afines [Monographs of
  the Institute of Mathematics and Related Sciences]}.
\newblock Instituto de Matem\'atica y Ciencias Afines, IMCA, Lima, 2003.

\bibitem{diff}
H.~Dulac.
\newblock Recherches sur les point singuliers des \'equations
  diff\'erentielles.
\newblock {\em J. Ecole Polytechnique}, (2):1--125, 1904.

\bibitem{Godement}
R.~Godement.
\newblock Th\'eorie des fais\c{c}eaux.
\newblock {\em Hermann, Paris}, 1973.

\bibitem{alcides}
Alcides Lins~Neto.
\newblock Construction of singular holomorphic vector fields and foliations in
  dimension two.
\newblock {\em J. Differential Geom.}, 26(1):1--31, 1987.

\bibitem{univ}
J.-F. Mattei.
\newblock Modules de feuilletages holomorphes singuliers. {I}.
  \'{E}quisingularit\'e.
\newblock {\em Invent. Math.}, 103(2):297--325, 1991.

\bibitem{MatQuasi}
J.~F. Mattei.
\newblock Quasi-homog\'en\'eit\'e et \'equir\'eductibilit\'e de feuilletages
  holomorphes en dimension deux.
\newblock {\em Ast\'erisque}, (261):xix, 253--276, 2000.
\newblock G\'eom\'etrie complexe et syst\`emes dynamiques (Orsay, 1995).

\bibitem{MM}
J.-F. Mattei and R.~Moussu.
\newblock Holonomie et int\'egrales premi\`eres.
\newblock {\em Ann. Sci. \'Ecole Norm. Sup. (4)}, 13(4):469--523, 1980.

\bibitem{MatSal}
J.-F. Mattei and E.~Salem.
\newblock Modules formels locaux de feuilletages holomorphes.
\newblock {\em PrePrint.}, 2004.

\bibitem{Seguy}
M.~Seguy.
\newblock {\em Cobordismes et reliabilit\'es \'equisinbuli\`eres de
  singularit\'es marqu\'ees de feuilletages holomorphes en dimension deux}.
\newblock PhD thesis, Toulouse, 2003.

\bibitem{reduction}
A.~Seidenberg.
\newblock Reduction of singularities of the differential equation
  {$A\,dy=B\,dx$}.
\newblock {\em Amer. J. Math.}, 90:248--269, 1968.

\bibitem{spivak}
Michael Spivak.
\newblock {\em A comprehensive introduction to differential geometry.}
\newblock Publish or Perish Inc., Boston, Mass., 1975.

\bibitem{zariski}
Oscar Zariski.
\newblock {\em Le probl\`eme des modules pour les branches planes}.
\newblock Hermann, Paris, second edition, 1986.
\newblock Course given at the Centre de Math\'ematiques de l'\'Ecole
  Polytechnique, Paris, October--November 1973.

\end{thebibliography}
\bibliographystyle{plain}
\end{document}